\documentstyle{amlts}
\begin{document}
\def\bye{\end{document}}
 \font\tenrm=cmr10

 \def\titleheadline#1{\def\one{#1}\ifx\one\empty\else
\gdef\thetitle{{\frenchspacing%
\let\\ \relax%\eightsc\uppercase{#1}
{#1}}}\fi}
\newif\ifshort
\def\shortname#1{\global\shorttrue\xdef
\theauthors{{\eightsc\uppercase{#1}}}}
\let\shorttitle\titleheadline

\def\joinrel{\mathrel{\mkern-4mu}}
\def\relbar{\mathrel{\smash-}}
\def\lrar{\relbar\joinrel\relbar\joinrel\relbar\joinrel\relbar\joinrel\rightarrow}
\def\srar{\relbar\joinrel\relbar\joinrel\rightarrow}

%\input BoxedEps.Tex %for mac
%\SetTexturesEPSFSpecial %for mac
\input boxedeps.tex % unix
\SetepsfEPSFSpecial % unix
\HideDisplacementBoxes
\def\figin#1#2{\medbreak
$$
 {\BoxedEPSF{#1.eps scaled
#2}%
}%
$$
\medbreak\noindent}
%--------------- Author macros ---------------
%for Bbb in amstex
\catcode`\@=11
\font\twelvemsb=msbm10 scaled 1100
\font\tenmsb=msbm10
%\font\ninemsb=msbm7 scaled 1100%msbm9
\font\ninemsb=msbm10 scaled 800
\newfam\msbfam
\textfont\msbfam=\twelvemsb  \scriptfont\msbfam=\ninemsb
  \scriptscriptfont\msbfam=\ninemsb
\def\msb@{\hexnumber@\msbfam}
\def\Bbb{\relax\ifmmode\let\next\Bbb@\else
 \def\next{\errmessage{Use \string\Bbb\space only in math
mode}}\fi\next}
\def\Bbb@#1{{\Bbb@@{#1}}}
\def\Bbb@@#1{\fam\msbfam#1}
\catcode`\@=12

 \catcode`\@=11
\font\twelveeuf=eufm10 scaled 1100
\font\teneuf=eufm10
\font\nineeuf=eufm7 scaled 1100%eufm9
\newfam\euffam
\textfont\euffam=\twelveeuf  \scriptfont\euffam=\teneuf
  \scriptscriptfont\euffam=\nineeuf
\def\euf@{\hexnumber@\euffam}
\def\frak{\relax\ifmmode\let\next\frak@\else
 \def\next{\errmessage{Use \string\frak\space only in math
mode}}\fi\next}
\def\frak@#1{{\frak@@{#1}}}
\def\frak@@#1{\fam\euffam#1}
\catcode`\@=12
 \annalsline{153}{2001}
\received{September 14, 1995}
\revised{November 9, 1999}
 \startingpage{297}

\font\emi= cmmi10 scaled 1700 
% Authors: Please start here:
%--------------- Author macros ---------------

%(Optional)
% Please enter all author-written macros 
% used in the body of the paper here:
% i.e.,
% \def\CC{\rm C\!\!\!\!I\,}, etc.

%-------------- Author entries --------------------

\title{On Brown-Peterson cohomology of \hbox{\emi QX}} %Article title

\shorttitle{ \eightsc\uppercase{On Brown-Peterson cohomology of} {\eightpoint \it QX}}

%% Please enter all acknowledgements here:
 \acknowledgements{Most of the work presented in the current paper was carried out during the
author's stay at  Kyoto University,  and partially supported by JSPS and
by a grant from the Ministry of Education of Japan.
}

% Please uncomment and use appropriate command:
\author{Takuji Kashiwabara}
%\twoauthors{}{}
%\authors{}% Separate each author with a comma and a space.

% Institution:
%% If more than one institution represented, please separate
%% with \\ , i.e.,
%% \institutions{University of Illinois at Chicago, Chicago, IL\\
%% Cornell University, Ithaca, NY}

\institutions{Institut Fourier,  Universit\'{e} de Grenoble  I, St.Martin-d'H\`{e}res,
France\\
{\eightpoint {\it E-mail address\/}: Takuji.Kashiwabara@ujf-grenoble.fr}}

%-------------- Article Text--------------------

%\intro %(Optional, Introduction)

\newcommand{\bpg}{\hat{\otimes }_{{\rm BP}^*}} 
\newcommand{\cotimes}{\hat{\otimes }}
\newcommand{\und}[2]{\underline{#1} _{#2}}
\newcommand{\Tor}{{\rm Tor}} \newcommand{\BP}[1]{{\rm BP}\langle{#1}\rangle}
\newcommand{\Hom}{{\rm Hom}} \newcommand{\colim}{{\rm colim}}
\renewcommand{\lim}{{\rm lim}}
\newcommand{\Ker}{{\rm Ker}} \renewcommand{\Im}{{\rm Im}}
\newcommand{\ra}{\rightarrow} \newcommand{\la}{\leftarrow}
\newcommand{\bea}{\begin{eqnarray*}} \newcommand{\eea}{\end{eqnarray*}}
\newcommand{\hra}{\hookrightarrow} \newcommand{\wt}{\widetilde}

\intro
 
Given a spectrum $X$ and a generalized cohomology theory $h$ with $h_*(X)$
known,  what can we say about $h_*(\underline{X} _i)$ where $\underline{X} _i$
denotes the $i^{\rm th}$ infinite loop space associated to $X$?  An obvious place to 
start investigating this question would be the case when $X$ is a suspension
spectrum $\Sigma ^{\infty }Y$ where $Y$ is a based space.  In this case one has
 $$\underline{X} _i\cong Q\Sigma ^iY\cong 
{\colim }_{s} \Sigma ^{i+s}\Omega ^sY.$$
The mod $p$ ordinary homology of such a space was computed by Kudo and Araki in the
case $p=2$ \cite{AK} and by Dyer and Lashof in the
case $p$ is odd \cite{dl}.  Later,  J. P. May  determined the Bockstein 
spectral sequence in terms of that for $Y$. Since the
 rational homology of such a space is easy to determine,
this gives us  complete knowledge of  ordinary homology of spaces of the form
$QX$.  As one might expect,  the first extraordinary homology that was studied
was the mod $p$ $K$-theory.
The first result is due to Hodgkin,
dating back to the  1970's.

\proclaimtitle{\cite{Ho}}
\proclaim{Theorem}
$K_*(QS^0;Z/p)\cong Z/p[\iota , Q\iota , Q^2\iota ,\dots ][\iota ^{-1}].$
\endproclaim

Here $Q^i$ denotes the $i^{\rm th}$ iteration of $Q$,  which is an analogue of
the classical Dyer-Lashof-Kudo-Araki operation,  defined up to a certain
indeterminacy.  Hodgkin had shown earlier in \cite{Ho2} that the indeterminacy was 
inevitable.  Later in the beginning of  the 1980's,  Miller and Snaith determined
  \cite{ms}
$K_*(QS^n;Z/2)$ as well as $K_*(QRP^n;Z/2)$.  Finally,  McClure succeeded in 
constructing a well-defined operation from $K_*(Y,Z/p^r)$ to $K_*(Y,Z/p^{r-1})$ 
($r\geq 2$),  for infinite loop spaces $Y$ and described $K_*(QX;Z/p^r)$ ($r\geq 1$)
in terms of the Bockstein spectral sequence for $K_*(X)$ (\cite{mc})
(which is equivalent to  the
knowledge of $K_*(X,Z/p^r)$ for all $r\geq 1$).  
The answer is too complicated to quote here,  but we note that for an
odd prime $p$, 
$K_*(QS^n;Z/p)$ is a free commutative (in the  graded sense) algebra with generators
$\iota , Q\iota ,Q^2\iota ,\dots $ where $\iota $ is the ``fundamental class",
that is, the image of the generator of $K_*(S^n)$ by the map induced by the 
map $S^n\ra QS^n$ with
$Q$ as above.

Naturally the next cases to 
study would be the higher Morava $K$-theories.  The computation by Hodgkin depended
on the theorem of Atiyah  that identifies the $K$-theory of
the classifying space of a finite 
group with the completion of its representation ring at the augmentation ideal
  $K^*(BG)\cong R(G)^{\wedge }$ and the intimate relationship between the 
classifying spaces of the symmetric groups and $QS^0$.  Unfortunately,  in the
case of Morava $K$-theories,  we only know that the  dimension of $K(n)^*(BG)$ is
equal to the number of $n$-tuples of commuting elements of order a power of $p$
\cite{HKR},
provided that it is concentrated in even degrees (which is the case for symmetric
groups (\cite{HKR}, \cite{Hu1})) and we do not have their functorial description.  
However, when $n=1$, it turned out
that these formulae were sufficient to recover the results of 
\cite{Ho}  (see \cite{K1}).  Furthermore,  by doing quite involved calculations
of characteristic classes,  the author computed $K(2)_*(QS^{r})$  (\cite{k2}).
The result is that one can define the operations $Q_1,\dots ,Q_{p+1}$ with 
indeterminacies,  which are subject only to the relations of the form
$Q_iQ_1=0$ if $1<i\leq p$,  $Q_{p+1}Q_1=Q_1Q_2$,  which have to be interpreted
suitably,  and $K(2)_*(QS^{r})$ is a free commutative algebra generated by 
iteration of these operations acting on the fundamental class.  

However,  the generalization of the method used in \cite{k2} seems to be 
difficult,  as there is a consensus among the experts that in the case for the
$n^{\rm th}$ Morava $K$-theory,  among the generalized Dyer-Lashof operations
there should be relations of length $n$, i.e., those of the form
$Q_{i_1}\cdots Q_{i_n}=\Sigma a_{i,j} Q_{j_1}\cdots Q_{j_n}$,  
that are not generated by shorter relations.  Since in the case for $n=2$,
the relations were derived from the computation of the image of $K(2)^*
(B\Sigma _{p^2})$ in $K(2)^*(B\Sigma _p\wr \Sigma _p)$,  a direct generalization
of the method in \cite{k2} would involve computing the image of
$K(n)^*(B\Sigma _{p^n})$ in $K(n)^*(\stackrel{n\mbox { -factors}}{\overbrace{B\Sigma _p\wr \cdots
\wr \Sigma _p})}$,  which seems to be 
 hopeless.

One thing that
one might want to try,  thus,  is to use $E_n$'s rather than $K(n)$'s,
so that one has generalized character theory (\cite{HKR}) at hand.  
Strickland and Turner  in \cite{ST} used this approach to describe 
$p^{-1}E_{n*}(CS^0)$ in terms of formal group laws. Here $CS^0$ is the disjoint
union of $B\Sigma _i$'s,  related to $QS^{2r}$'s by Thom isomorphisms and 
Snaith splitting,     whose group completion has the homotopy type of $QS^0$.
  Furthermore,  Strickland carried out more detailed analysis of
$E_{n}^*(CS^0)$ using rich structures of $CS^0$  to give its formal group 
theoretic description \cite{S}.  After the first version of this paper was 
submitted,  Strickland went on to obtain (\cite{neil}) a functorial description of
$$\pi_*(L_{K(n)}(E_n\wedge QX))$$ when $\pi_*(L_{K(n)}(E_n\wedge X))$ is the
completion of a free module concentrated in even degrees.  Here  
$\pi_*(L_{K(n)}(E_n\wedge -))$ replaces $E_{n*}(-)$.  

On the other hand,  the result by Kudo and Araki can be reinterpreted in terms of
cohomology using a  result due to Lannes and Zarati \cite{LZ}. 
Let $D$ be the destabilization functor, i.e., the left adjoint to
the forgetful  functor from the category of unstable algebras
 over the Steenrod algebra
 to the category of modules over the Steenrod algebra, and $D^i$ its $i^{\rm th}$  left derived
functor.  Then their results imply the isomorphism
 $PH^*(QX;Z/2)\cong \oplus _i\Sigma ^2D^i\Sigma ^ {2-i}\tilde{H} ^* (X)$ where 
$P$ denotes the module of primitives, although they do not state it in this way. 
This isomorphism can be seen using the relationship between 
$PH^*(QX;Z/2)$ and $R_nH^*(X;Z/2)$ 
which
seems to be well-known to experts but whose proof does not appear in the 
literature, where $R_n$ is as defined in \cite{LZ}.
We will discuss it in full detail in a
subsequent work \cite{Kprep}. 

 This result suggests that in the case of 
generalized cohomology theories too,  cohomology operations should play a major 
role.  In \cite{bpqs} the author introduced the notion of 
the destabilization functor for the ${\rm BP}$-cohomology.  It is not a precise 
analogue of $D$ above,  but an analogue of the composition $U\circ D$ where $U$
is the universal enveloping unstable algebra functor,  that is, the left 
adjoint to the forgetful functor from the category of unstable algebras over
the Steenrod algebra to that of unstable modules over the Steenrod algebra.
More precisely:

\numbereddemo{Definition}
Denote by $D_{\rm BP}$ the left adjoint to the
forgetful  functor from the category of ${\rm BP}$-unstable algebras introduced by
Boardman, Johnson, and Wilson (\cite{bjw}),  to the category of  modules over
the Landweber-Novikov algebra,  ${\rm BP}^*{\rm (BP)}$.
\enddemo

And the author proved:
 
\proclaimtitle{\cite{bpqs}} 
\proclaim{Theorem}
\label{old;main}
Let $X$ be a $(-1)$\/{\rm -}\/connected
spectrum whose stable cells are concentrated in even degrees and 
whose ${\rm BP}$\/{\rm -}\/cohomology is finitely generated as a ${\rm BP}^*{\rm (BP)}$\/{\rm -}\/module{\rm .}
  Then the
natural map $D_{\rm BP}({\rm BP}^*(X))\ra {\rm BP}^*(\Omega ^{\infty }X)$ is an isomorphism of
${\rm BP}^*$\/{\rm -}\/algebras{\rm .}
\endproclaim 

The statement and the proof in \cite{bpqs} contain minor errors;  see Section 7.
This result notably applies
to ${\rm BP}^*(QS^{2j})$.  We also note that there are many other infinite loop spaces
that fit into this category   whose mod $p$ ordinary cohomology is still 
unknown. Similar results concerning $K$-theory were also proved by Bousfield,
namely:
\proclaimtitle{\cite[Th.\ 8.3, Cor.\ 8.6]{BK}}
\proclaim{Theorem}
If $E$ is a spectrum with $K^*(E;Z^{\wedge }_p)$ torsion\/{\rm -}\/free{\rm ,} then 
$K^*(\Omega ^{\infty }_0E;Z^{\wedge }_p)$ is naturally isomorphic to $WK^*(E;Z^{\wedge }_p)_H${\rm .}
Furthermore{\rm ,} if $E$ is $0$\/{\rm -}\/connected and $H^i(E;Z^{\wedge }_p)=0 $ for $i=1,2${\rm ,}
then\break $K^*(\Omega ^{\infty }E;Z^{\wedge }_p)$ is naturally isomorphic to $TK^*(E;Z^{\wedge }_p)${\rm .}
Here $T$ is the free\break $\theta ^p$\/{\rm -}\/ring functor{\rm ,} i{\rm .}e{\rm .,} the left adjoint to the
forgetful functor from the category of $\theta ^p$\/{\rm -}\/rings to that of
profinite $p$\/{\rm -}\/groups{\rm ,}  $K^*(E;Z^{\wedge }_p)_H$ is a certain enrichment of
$K^*(E;Z^{\wedge }_p)${\rm ,}  and $W$ is an appropriate adjoint functor that takes
into account this enrichment{\rm .}
\endproclaim

We will discuss the relationship between the Bousefield functors $T$, $W$ and the
$K$-theory version of our destabilization functor in \cite{Kprep}. 
Furthermore,   we know the algebra
structure of $K(n)_*(QS^i)$'s:

\proclaimtitle{\cite{bpqs}}
\proclaim{Theorem}\label{old;free} $K(n)_*(QS^{2m})$ $(m\geq 0)$
is a polynomial algebra concentrated in even degrees{\rm .}  $K(n)_*(QS^0)$ is a 
tensor product of a polynomial algebra concentrated in even degrees with 
$K(n)_*[Z]${\rm .}  $K(n)_*(QS^{2m-1})$ is an exterior algebra with generators in odd
degrees{\rm .}  \endproclaim

It was   shown in \cite{bcat} that these algebras are also cofree as coalgebras.  On 
the other hand, their Hopf algebra structures still remain to be studied.

Now,  the purpose of this paper is to generalize the results above.  
First of all, as we deal with the spaces that are not necessarily finite,
we should take into account the topology on their ${\rm BP}$-cohomology.  
Unfortunately abelian toplogical groups do not form an abelian category, so we
need to take care to set up a correct framework to define the destabilization
functor.  

\numbereddemo{Definition}
We will call ${\cal M}_{\rm BP}$ and ${\cal K}_{\rm BP}$ respectively the
category of stable ${\rm BP}$-cohomology modules defined in \cite{boa} and that
of unstable ${\rm BP}$-cohomology algebras defined in \cite{bjw} respectively, with  
the following additional requirements:
\medbreak
\listitem{i} The filtration is over $Z$.
\smallbreak
\listitem{ii} The elements of degree $n$ have  filtration at least $n$.
\medbreak
Now, denote   ${\cal K}_{_0 {\rm BP}}$ the category of augmented unstable 
${\rm BP}$-cohomology algebras,  that is, the category whose objects are 
 unstable 
${\rm BP}$-cohomology algebras equipped with the augmentation to the coefficient
ring ${\rm BP}^*$, and whose morphisms are morphisms of unstable 
${\rm BP}$-cohomology algebras that respect the augmentation.
\enddemo
 
The Brown-Peterson cohomology of a spectrum or a pointed space (not necessarily
of finite type) is equipped with
the skeletal filtration and  become  objects of ${\cal M}_{\rm BP}$ or ${\cal K}_{_0 {\rm BP}}$
respectively.  (The systematic use of the skeletal filtration causes some 
inconvenience.  However,  it allows us to deal with the  categorical sums more
easily.)
Thus ${\cal M}_{\rm BP}$ is the algebraic model for the category of spectra,
${\cal K}_{_0 {\rm BP}}$ that of pointed spaces,  and
the augmentation ideal functor
${\cal I}$  from ${\cal K}_{_0{\rm BP}}$ to 
${\cal M}_{\rm BP}$ corresponds to the suspension spectrum functor
 $\Sigma ^{\infty }$.  Note also that the cokernel exists in the categories 
${\cal M}_{\rm BP}$ and ${\cal K}_{_0{\rm BP}}$.  It is nothing but the algebraic cokernel
equipped with the quotient filtration.

\numbereddemo{Definition}
A stable ${\rm BP}$-cohomology module is called {\it free} if it is in the essential image
of the left adjoint to the forgetful functor to the category of graded sets.
A stable ${\rm BP}$-cohomology module is called {\it well\/{\rm -}\/presented} if it is isomorphic in
${\cal M}_{\rm BP}$ to the cokernel of a map of the form $F^1\ra F^0$ where $F^i$'s
are free.  We name the full subcategory of 
${\cal M}_{\rm BP}$ formed by well-presented objects ${\cal M'}_{\rm BP}$.
An unstable ${\rm BP}$-cohomology algebra is said to be well-presented if
it is well-presented as a stable ${\rm BP}$-cohomology module.  We call the full
subcategory of ${\cal K}_{_0{\rm BP}}$ formed by well-presented objects ${\cal K'}_{
_0{\rm BP}}$.
\enddemo

The point is that in general, as far as the module structure is concerned,
 one can express any module as a quotient of a
free one, but there is no guarantee that the original topology would coincide 
with the quotient topology.  We avoid this problem by considering only those
modules that have the ``correct'' topology.
Note also that if $X$ is a finite-type wedge of suspensions of ${\rm BP}$,  then
${\rm BP}^*(X)$ is free.  If $X$ is an arbitrary wedge of  suspensions of ${\rm BP}$,  
i.e., $X\cong \vee _i \Sigma ^{d_i}{\rm BP}$,  then the completion with respect to the
skeletal topology of the direct sum $\oplus _i {\rm BP}^*(\Sigma ^{d_i}{\rm BP})$ is free,
and it is contained in ${\rm BP}^*(X)$.  It is dense if one considers the 
finite-subcomplex topology defined in Section \ref{bp;des}.  
Now we can introduce our destabilization functor:

\proclaim{Proposition-Definition}
The augmentation ideal functor from ${\cal K'}_{_0{\rm BP}}$ to ${\cal M'}_{\rm BP}$ admits a left
adjoint{\rm ,} which we call the destabilization functor and we note it ${\cal D}${\rm .}
\endproclaim

For free objects in ${\cal M'}_{\rm BP}$ that are
 completions of $\oplus _i {\rm BP}^*(\Sigma ^{d_i}{\rm BP})$,  one defines the value of
${\cal D}$ to be the completions of $\otimes _i {\rm BP}^*(\und{{\rm BP}}{d_i})$, where
$\und{{\rm BP}}{d_i}\cong \Omega ^{\infty }\Sigma ^{d_i}{\rm BP}$.  Note that ${\cal D}$ is isomorphic
to ${\rm BP}^*(\Pi _i(\und{{\rm BP}}{d_i}))$ if $\Pi _i(\und{{\rm BP}}{d_i})$ (equivalently 
$\vee _i \Sigma ^{d_i}{\rm BP}$) is of finite type.  Otherwise it is dense in
${\rm BP}^*(\Pi _i(\und{{\rm BP}}{d_i}))$ with respect to the
finite-subcomplex topology.  
Since 
any object in ${\cal M'}_{\rm BP}$ is a quotient of a free one,
  we can define ${\cal D}$ noting that it has to be right exact.
Of course, we will need to know if the ${\rm BP}$-cohomology of our spaces are 
well-presented.  For this purpose, we introduce yet another definition that is
easier to check:

\numbereddemo{Definition}
Let $X$ be a space or a spectrum.  ${\rm BP}^*(X)$ is said to be {\it well\/{\rm -}\/generated} if 
${\rm BP}^*(X)\bpg (Z/p)\hra H^*(X;Z/p)$.  \enddemo

In most cases when the ${\rm BP}$ cohomology of a space is known, it is well-generated.
However, the result of Inoue \cite{In} seems to imply that
${\rm BP}^*(B{\rm SO}(6))$ is not.  
Anyhow,  we will prove the following:

\proclaim{Lemma} \label{gen=pre}
For a spectrum or space
$X${\rm ,} if ${\rm BP}^*(X)$ is well\/{\rm -}\/generated{\rm ,} then it is well\/{\rm -}\/presented{\rm .}
\endproclaim

Now we are ready to state the generalization of the preceding results.
It is quite natural to limit ourselves to the spaces $X$ such that
${\rm BP}^*(X)$ satisfies   Landweber's exact-functor-theorem's
hypothesis to have   Kunneth's isomorphism at hand.  Thus our
 first main result is that the functor $Q$ preserves these properties.  That is:

\proclaim{Theorem} \label{qx;ten}
Let $X$ be a connected space such that
${\rm BP}^*(X)$ is\break Landweber\/{\rm -}\/flat{\rm ,} i.e.{\rm ,} it satisfies one of the equivalent conditions in Theorem
{\rm \ref{land;gen}.}  Then ${\rm BP}^*(QX)$ satisfies the same conditions{\rm .}  Furthermore{\rm ,}
if ${\rm BP}^*(X)$ is well\/{\rm -}\/generated{\rm ,} then so is ${\rm BP}^*(QX)${\rm .}
\endproclaim

Now that we are assured that our ${\rm BP}^*(QX)$ lies in the right category, we can
compare it with ${\cal D}\wt{{\rm BP}}^*(X)$ and we have:

\proclaim{Theorem} \label{realmain}
Let $X$ be a connected space{\rm ,}  satisfying one of the equivalent conditions in Theorem
{\rm \ref{land;gen},}   whose ${\rm BP}$ cohomology is well\/{\rm -}\/generated{\rm .}
Then the natural map ${\cal D}\wt{{\rm BP}}^*(X)\ra {\rm BP}^*(QX)$ is an isomorphism
in ${\cal K'}_{_0{\rm BP}}${\rm .}
\endproclaim

We can also generalize Theorem \ref{old;main}:

\proclaim{Theorem} \label{evencell}
Let $X$ be a  {\rm ($-1$)-}\/connected spectrum which has stable cells only in even
degrees{\rm .}  Then the natural map ${\cal D}{\rm BP}^*(\Sigma ^iX)\ra {\rm BP}^*(\und{X}{i})$
is an isomorphism if $i\geq 0${\rm .}
\endproclaim

 What this
means is as follows:  $Q$ factors as the composition $\Omega  ^{\infty}
\Sigma ^{\infty}$ where $\Sigma ^{\infty}$ is the functor which
associates to a space its suspension spectrum,  $\Omega  ^{\infty}$ its right
adjoint.  On the other hand,  ${\rm BP}$-cohomology of a spectrum takes value in the
category in the modules over ${\rm BP}^*{\rm (BP)}$ whereas ${\rm BP}$-cohomology of a space
takes value in the category of ${\rm BP}$-unstable algebras,  and $\Sigma ^{\infty}$ 
is  compatible with the augmentation ideal functor ${\cal I}$,  which means that ${\cal D}$
is an algebraic model for $\Omega ^{\infty }$.  Thus the composition
${\cal D}\circ {\cal I}$ can be regarded as 
an algebraic counterpart of the composition of the functors $X\mapsto 
\Sigma ^{\infty}X\mapsto QX$.  Our theorem shows that it is a good model.  This
description of ${\cal D}$ looks quite abstract.  However,  it can be described
completely algebraically and concretely.  As a matter of fact we first obtain
an algebraic answer for Morava $K$-theories and ${\rm BP}$-cohomology as follows,
then identify our answer with the result of the
destabilization.

\proclaim{Theorem} \label{knqx;alg}
Let $X$ be a space as in Theorem {\rm \ref{realmain}.}
Furthermore let $\{ f_i:X\ra \und{\rm BP}{d_i}|i\in I\}$ be a set of 
topological ${\rm BP}^*{\rm (BP)}$\/{\rm -}\/module
generators for $\wt{{\rm BP}^*}(X)$ {\rm (}\/that is{\rm ,} the ${\rm BP}^*{\rm (BP)}$\/{\rm -}\/submodule
generated by $f_i$\/{\rm '}\/s is dense\/{\rm ),}\break   $\{ g_j:\vee _i\Sigma ^{d_i}BP\ra
\vee _j\Sigma ^{e_j}BP
\} $ be a set that generates topologically
 a complete set of relations {\rm (}i.e{\rm .,}    the 
exact sequence $0\la \wt{{\rm BP}^*}(X)\la {\rm BP}^*(\vee _i\Sigma ^{d_i}{\rm BP})\la
{\rm BP}^*(\vee _j\Sigma ^{e_j}{\rm BP})$ exists{\rm ).}  Then there are
{\rm \listitem{i}} an exact sequence of Hopf algebras
$$K(n)_*\ra K(n)_*(QX) \ra K(n)_*(\Pi _i\und{\rm BP}{d_i})\ra
K(n)_*(\Pi _j\und{\rm BP}{e_j})$$
{\rm \listitem{\rm ii}} and a coexact sequence of algebras
$${\rm BP}^*\la {\rm BP}^*(QX) \la {\rm BP}^* (\Pi _i\und{\rm BP}{d_i})\la
{\rm BP}^*(\Pi _j\und{\rm BP}{e_j}).$$
\endproclaim

Our description of $K(n)_*(QX)$ and ${\rm BP}^*(QX)$
in Theorem \ref{knqx;alg} may not look algebraic.  However, we will explain how
we can reduce everything to a pure algebra (provided that one has the 
 ${\rm BP}^*{\rm (BP)}$-module presentation of $\wt{\rm BP} ^*(X)$,  but this is again a
purely algebraic question) using the determination of $E_*(\und{\rm BP}{*})$ for
complex oriented homology theories $E$ by Ravenel and Wilson \cite{RW}.  As far
as the algebra and coalgebra structure is concerned,  one can be more explicit
(again, the following will be proved before being used in the proof of Theorem
0.14) and generalize Theorem \ref{old;free}.

\proclaim{Theorem} \label{int:knfree}
Let $X$ be as above{\rm .}  Then $K(n)_*(QX)$ is a free commutative algebra{\rm .}
Furthermore{\rm ,} $K(n)_*(Q\Sigma X)$ is a cofree cocommutative coalgebra{\rm .}
\endproclaim

Another natural question to ask here is how the set of
${\rm BP}^*$-module generators of
${\rm BP}^*(QX)$ for such spaces (since we know by Theorem \ref{qx;ten} that they get detected
by mod $p$ ordinary cohomology) can be described in terms of $H^*(QX)$. 
  Denote by $M_X$ the image of the Thom map
$\rho _X: {\rm BP}^*(X) \ra H^*(X)$.    Then this is equivalent to the knowledge of
$M_{QX}$.
Our  answer to this question is as follows:

\proclaim{Theorem} \label{int:qbpg}
Let $X$ be a connected space with the property that ${\rm BP}^*(X^j)\cong {\rm BP}^*(X)^{\hat{\otimes }j}${\rm ,}
and ${\rm BP}^*(X)\bpg Z/p\hra H^*(X;Z/p)${\rm .}
Then $M_{QX}\break =C$ where $C$ will be defined as in Proposition
{\rm \ref{bpg;det}.}
\endproclaim

In this paper we use the following convention.  ${\rm BP}$ will denote the\break 
$p$-complete version of the Brown-Peterson spectrum (what would normally be\break denoted
as ${\rm BP} _p^{\wedge }$) with the coefficient ring ${\rm BP}_*\cong Z^{\wedge }_p[v_1,
v_2, \dots ,]$,  $|v_n| =\break 2(p^n-1)$,  $K(n)$ the usual $n^{\rm th}$ Morava $K$-theory 
with $K(n)_*=Z/p[v_n,v_n^{-1}]$,  $E(n)$ the Johnson-Wilson theory
with $E(n)_*=Z_{(p)}[v_1,\dots ,v_n,v_n^{-1}]$,  $H$ the mod $p$ ordinary  
(co)homology.  Throughout the main text of the paper,  $p$ will be an odd prime.  However,
most of our results also hold for $p=2$.  Necessary modifications are indicated
in the appendix.  A ``space" will mean a pointed  topological
space with the
homotopy type of a CW-complex of finite type unless otherwise specified.  
A generalized cohomology of
a space is topologized via the skeletal filtration unless otherwise specified,
and it is with respect to this topology that we take the completed tensor 
products.  Thus we ignore the ``$p$-adic part" of the inverse limit topology on
${\rm BP}^*(X)\cong\lim_{n,i}{\rm BP}^*(sk_nX)/p^i$.  As a matter of fact, this does not
make much difference since we take into account the module structure over 
${\rm BP}^*$, thus over the $p$-adics, and the finite type hypotheses imply that 
the skeleton topology together with the module structure over the $p$-adics
suffices to determine the inverse limit topology.

The organization of this paper is as follows.  In Section \ref{prel} we collect
the facts on $QX$ that are necessary for us.
In Section \ref{fil;dll},  we define the notion of Dyer-Lashof length-like
filtration,  which is used repeatedly in the paper.  
In Section \ref{rwy;rev} we review and generalize relevant
results in \cite{rwy}.  In Section \ref{wreath} we use
the results in Section \ref{rwy;rev} and a result by Hunton on the behavior of
the Atiyah-Serre-Hirzebruch spectral sequence for a wreath product to show that
many properties that ${\rm BP}^*(X)$ possesses are passed onto ${\rm BP}^*(D_pX)$ and thus to
${\rm BP}^*(QX)$.  In Section \ref{thom;im} 
under the assumption that ${\rm BP}^*(X)\bpg (Z/p)\subset H^*(X)$ and that 
${\rm BP}^*(X)$ satisfies   Landweber's exact-functor-theorem's hypothesis,
we determine the image of the Thom map ${\rm BP}^*(QX)\rightarrow H^*(QX)$.  
In Section \ref{kn;poly} we use these results to conclude 
that $K(n)_*(QX)$ injects to a product of $K(n)_*(\underline{\rm BP} _i)$'s,  
and deduce from it 
 that $K(n)_*(QX)$ is a free commutative 
algebra.  Then we proceed further
to show that the cokernel of the map $K(n)_*(QX)
\rightarrow \otimes K(n)_*(\underline{\rm BP} _i)$ again injects to a
product of $K(n)_*(\underline{\rm BP} _i)$'s, and get a completely algebraic 
description of these objects.

The author would like to
thank Steve Wilson,  and John Hunton as well as many other people
for helpful conversation.  The author also thanks the referee for 
careful proofreading and suggestions on the exposition.

\section{Preliminaries} \label{prel}

In this section we collect mostly well-known facts on infinite loop spaces needed later in the paper.

\numbereddemo{Definition} \label{def;11}
Let $I=(\varepsilon _1,s_1,\dots \varepsilon _k,s_k)$ such that $s_j\geq \varepsilon _j$ 
and $\varepsilon _j= 0$ or $1$.  Define the degree (d),  the excess (e),  
the length  (l), and the presence of Bockstein at the end (b) of $I$ by
\bea
\noalign{\vskip4pt}
d(I) & = & \Sigma _{j=1}^k[2(s_j(p-1)-\varepsilon _j] ,\\ \noalign{\vskip4pt}
e(I) & = & 2s_1-\varepsilon _1-\Sigma _{j=2}^k[2(s_j(p-1)-\varepsilon _j],\\ \noalign{\vskip4pt}
l(I) & = & k,\\ \noalign{\vskip4pt}
b(I) &=&\varepsilon _1 .\\
\noalign{\vskip-9pt}
\eea
$I$ is said to be {\it admissible} if $ps_j-\varepsilon _j\geq s_{j-1} $ for $2\leq j\leq k$.  
\enddemo

For any such sequence $I$ (not necessarily admissible),  one has a corresponding
homology operation on $E_{\infty }$ spaces $Q^I=
\beta ^{\varepsilon _1}Q^{s_1}\dots \dots  \beta ^{\varepsilon _k}Q^{s_k}$,
that raises the degree of elements by $d(I)$ and vanishes on elements of
degree greater than $e(I)$.

\proclaimtitle{\cite{dl}}
\proclaim{Theorem}
 \label{th;12}
Let $X$ be a connected space{\rm ,}  and $\Lambda $ a basis for $\widetilde{H} _*(X)${\rm .}
Then
$H_*(QX)$ is a free commutative algebra on generators
$Q^I(x)$  {\rm (}$x\in \Lambda ${\rm ),}  where $I$ is admissible{\rm ,}  $e(I)+b(I)>{\rm deg}(x)${\rm .}  Now{\rm ,}
$H_*(CS^0)$ is  a free commutative algebra on generators
$Q^I([1])${\rm ,}   where $I$ is admissible{\rm ,}  $e(I)>0${\rm .}
Finally{\rm ,}  $H_*(QS^0)$ is the algebra generated by
$H_*(CS^0)${\rm ,}  and $[-1]${\rm ,}  subject to the relation $[1]\cdot [-1]=1${\rm .}
In the above{\rm ,}
 $[i]$ denotes
the image of the element $i\in \pi _0(QS^0)\cong Z$ or the element $i\in \pi _0(
CS^0)\cong Z^+$  by
Hurewicz homomorphism $\pi _0(-)\rightarrow H_0(-)${\rm .}
\endproclaim

Since the sphere spectrum is a ring spectrum,  its multiplication induces a
pairing $QS^{i}\times QS^{j}\rightarrow QS^{i+j}$.  When  $i=j=0$,
this pairing agrees with the map induced by $CS^0\times CS^0\rightarrow CS^0$
whose components are given by the maps induced by $\Sigma _a\times \Sigma _b
\rightarrow \Sigma _{ab}$.  This induces a pairing in homology,  denoted by 
$\circ $: $H_*(CS^0)\otimes H_*(CS^0)\rightarrow H_*(CS^0)$.  
Since any spectrum,  thus in particular a suspension spectrum,  is a module
spectrum over the sphere spectrum,  one gets a pairing $QS^0\otimes QX\ra
QX$,  which often is called the composition pairing,  as it agrees with the
colimit of the maps given by the composition 
\begin{eqnarray*}
\noalign{\vskip4pt}
\Omega ^nS^n\times \Omega^n \Sigma ^nX&=&{\rm Map} _*(S^n,S^n)\times 
{\rm Map} _*(S^n,\Sigma ^nX)\\
\noalign{\vskip4pt}
&& \ra {\rm Map} _*(S^n,\Sigma ^nX) =\Omega^n \Sigma ^nX
. \\
\noalign{\vskip-9pt}
\end{eqnarray*}
  We still denote the induced pairing in homology 
$H_*(QS^0)\otimes H_*(QX)\rightarrow H_*(QX)$ by $\circ $.
Furthermore,
the usual Pontrjagin product of $H_*(QX)$ will be denoted by $\star $
or just by juxtaposition.  We will need the following. 

\proclaimtitle{cf.\ \cite{May}} 
\proclaim{Theorem}
\label{mayf}
\medbreak
{\rm \listitem{i}} {\rm (}\/May\/{\rm '}\/s formula\/{\rm )} $Q^s[1]\circ x=\Sigma _{t\geq 0}Q^{s+t}P^t_*x$
\medbreak
{\rm \listitem{ii}} {\rm (}\/May\/{\rm '}\/s formula\/{\rm )} $\beta Q^s[1]\circ x=\Sigma _{t\geq 0}\beta
Q^{s+t}P^t_*x
	-\Sigma _{t\geq 0}Q^{s+t}P^t_*\beta x$
\medbreak{\rm \listitem{iii}} {\rm (}\/Nishida relation\/{\rm )} $P^r_*(Q^s(x))=\Sigma _i(-1)^{r+i} {(p-1)(s-r) \choose
r-pi}
	 Q^{s-r+i}P^i_*(x)$
\medbreak
{\rm \listitem{iv}}\/ {\rm (}\/Nishida relation\/{\rm )} \bea P^r_*\beta (Q^s(x))&=&
	\Sigma _i(-1)^{r+i}{(p-1)(s-r)-1 \choose r-pi}\beta Q^{s-r+i}P^i_*(x)
	\\ \noalign{\vskip4pt} & & + \ \Sigma _i(-1)^{r+i}{(p-1)(s-r)-1 \choose r-pi-1}
		Q^{s-r+i}P^i_*\beta (x)\eea
\medbreak{\rm \listitem{v}} $Q^s\sigma(x)=\sigma( Q^s(x))$
\medbreak{\rm \listitem{vi}} $Q^s(x)=x^p $ if $s=(p-1){\rm deg}(x)${\rm .}
\medbreak
\noindent Here $P^r_*$ denotes the dual Steenrod reduced power operation{\rm .}
\endproclaim

Note that the formulae i) and ii) can be combined into a single formula:
$$\beta ^{\varepsilon }Q^s[1]\circ x=\Sigma _{t\geq 0}\beta^{\varepsilon } Q^{s+t}
P^t_*x
        -\varepsilon \Sigma _{t\geq 0}Q^{s+t}P^t_*\beta x.$$
Next we need to know $H^*(B\Sigma _p)$.

\proclaimtitle{e.g. \cite{SE}}
\proclaim{Proposition}
 \label{prop;14}
$H^*(B\Sigma _p)\hra H^*(BZ/p)=\Lambda (x)\otimes Z/p[y]${\rm ,}  the image is the
subalgebra generated by $y^{p-1}$ and $xy^{p-2}${\rm .}  
\endproclaim

We denote by $e_{2i(p-1)}$ the element in $H_*(B\Sigma _p)$ or $H_*(BZ/p)$
that is dual to $y^{i(p-1)}$.   Then by 
definition \cite{AK}, \cite{dl},  in $H_*(B\Sigma _p)\subset H_*(CS^0)$,  
$Q^i[1]=e_{2i(p-1)}$.  We note that,  since $\beta (x)=y$,  an easy 
Atiyah-Hirzebruch spectral sequence argument shows that the image of 
${\rm BP}^*(B\Sigma _p)\ra H^*(B\Sigma _p)$ is the subalgebra generated by $y^{p-1}$.

Finally,  we will need the following:

\numbereddemo{Definition} \label{combdef} {\rm
For a topological space $X$,  $CX=\mathbold{\cup}_nE\Sigma _n\times _{\Sigma _n}X^n/\sim$, where $\sim$ is the
equivalence defined in \cite{geo}.  Define a filtration $F$  on
$CX$ by $$F_i(CX)=\mathbold{\cup}_{n\leq i}\Im (E\Sigma _n\times _{\Sigma _n}X^n\rightarrow
CX).$$  Define 
$DX$ to be  $\vee _i F_i(CX)/F_{i-1}(CX)$.}
\enddemo

\demo{Remark} It is well-known and easy to see that $ F_i(CX)/F_{i-1}(CX)$ 
is homeomorphic to $E\Sigma _n^+ \wedge X^n$.\enddemo

\proclaimtitle{e.g.\ \cite{sn}, \cite{be}} 
\proclaim{Theorem}
 For a connected space $X${\rm , \label{snsp}}
$QX$ is stably homotopy equivalent to $DX${\rm ,}  and $CX$ is homotopy equivalent to
$QX${\rm .}
\endproclaim

\section{Dyer-Lashof length-like filtration} \label{fil;dll}

In this section,  we introduce the notion of Dyer-Lashof length-like filtration,
and discuss its properties.

\numbereddemo{Definition}
Let $A$ be an algebra augmented over a field $k$.  By ${\rm Ind}(A)$ we denote
the {\it module of indecomposables}\/,  i.e., ${\rm Ind}(A)=I/I^2$,  where $I=\Ker (A\ra k)$.
\enddemo

We use this notation instead of the traditional $Q$ to avoid confusion with
other $Q$'s used   throughout the paper.

\numbereddemo{Definition} \label{filt;dll}
Fix $d\in Z$.  Denote by $Z_{\geq d}$ the set of integers greater than or equal
to $d$.  Let $\{ A_{i,*}|i\in Z_{\geq d}\}$ be a family of 
bi-graded Hopf algebras
augmented over a field $k$ with characteristic $p$.  
Suppose that this family is equipped with a suspension map,
i.e.,  a morphism of $k$-vector spaces $\sigma: A_{i,*}\ra A_{i+1,*+1}$.  
We say that the family $\{ A_{i,*}|i\in Z_{\geq d}\}$ together with $\sigma $
form a $Z_{\geq d}$-indexed family of graded algebras with Dyer-Lashof length-like filtration if they are equipped with an
increasing filtration $F$ on  each $A_{i,*}$,
$k=F_0(A)\subset F_1(A)\subset \cdots \subset F_n(A)\subset \cdots \subset A$
satisfying the following properties:
\medbreak
\listitem{i} Each $A_{i,*}$ is a free commutative algebra (in a graded sense).
\smallbreak\listitem{ii} $\sigma $ factors through ${\rm Ind}(A)$, and its image is contained in $PA$
.
\smallbreak\listitem{iii} The decomposition $A_{i,*}\cong P_{i,*}\otimes E_{i,*}$ with $P$ polynomial and $E$ 
exterior holds as Hopf algebras.  
\smallbreak\listitem{iv}$\sigma $ induces an isomorphism ${\rm Ind}P_{i,*}\ra 
{\rm Ind}E_{i+1,*+1}$, an injection ${\rm Ind}E_{i,*}\ra\break P_{i+1,*+1}$.  
\smallbreak\listitem{v} $F$ is exhaustive;  i.e., $\mathbold{\cup}_jF_j(A_{i,*})=A_{i,*}$.
\smallbreak\listitem{vi} $F_l(A_{i})\cdot F_k(A_i)\subset F_{l+k}(A_i)$.
\smallbreak\listitem{vii} $A_i$ is isomorphic as an algebra
to its associate graded object with respect to the 
filtration $F$ with the induced multiplication.
\smallbreak\listitem{viii} $F$ is compatible with $\sigma $;  i.e., $\sigma (F_j(A_i))\subset
F_{j}(A_{i+1})$.
\smallbreak\listitem{ix} $\sigma $ induces isomorphisms $$F_1(A_{i,*})/F_0(A_{i,*})\cong
\colim (A_{i,*}\stackrel{\sigma }{\ra } 
A_{i+1,*+1}\stackrel{\sigma }{\ra }\cdots ).$$
\enddemo

\demo{Remark} \label{dll;name}
Let $A_i$ be the subalgebra of $H_*(QS^i;Z/p)$ generated by the elements of the
form $Q^I(\iota _i)$ where $\iota _i$ is the fundamental class in $H_i(QS^i;Z/p)$
and where $I$ contains no Bockstein.  Then it becomes a  $Z^+$-indexed family of
algebras with Dyer-Lashof length-like filtration by defining $\sigma $ to be the
restriction of the homology suspension map and $F$ 
by $F_m(A)$ to be the span of monomials of weight less than or equal to $m$,
where the weight is defined by ${\rm weight}(Q^I(\iota ))=p^{l(I)}$,  ${\rm weight}(xy)
={\rm weight}(x)+{\rm weight}(y)$.  This is the origin of the name of Dyer-Lashof length-like filtration.\enddemo

The following three propositions are straightforward consequences of the definition.

\proclaim{Proposition} \label{dll;ten}
Let $A${\rm ,} $B$ be $Z_{\geq d}$\/{\rm -}\/indexed families of algebras equipped with a 
Dyer\/{\rm -}\/Lashof length-like filtration{\rm .} Then so is $A\otimes B$ with the suspension
given by $$\sigma:A\otimes B\ra {\rm Ind}(A\otimes B)\cong {\rm Ind}(A)\oplus {\rm Ind}(B)
\stackrel{\sigma \oplus \sigma }{\ra }A\oplus B \ra A\otimes B $$  and the usual tensor product
filtration.
\endproclaim

\proclaim{Proposition} \label{dll;col}
A direct limit of inclusions  compatible with the filtration and the
suspension of $Z_{\geq d}$\/{\rm -}\/indexed families of algebras equipped with a 
Dyer\/{\rm -}\/Lashof length-like filtration is again a $Z_{\geq d}$\/{\rm -}\/indexed family
of algebras equipped with a Dyer\/{\rm -}\/Lashof length-like filtration{\rm .}
\endproclaim

\proclaim{Proposition}\label{dll;quo} Let $A$ be a $Z_{\geq d}$\/{\rm -}\/indexed family of algebras equipped with 
a Dyer\/{\rm -}\/Lashof length-like filtration{\rm .} 
If $I_i$\/{\rm '}\/s form a family of ideals of $A_i$\/{\rm '}\/s compatible with the suspension such
that $A_i/I_i$\/{\rm '}\/s are free commutative algebras{\rm ,}  then $A_i/I_i$\/{\rm '}\/s form a
$Z_{\geq d}$\/{\rm -}\/indexed family of algebras equipped with a
 Dyer\/{\rm -}\/Lashof length-like filtration{\rm .}
\endproclaim

For future use,  we record some examples.

\proclaim{Proposition} \label{dll;qx}
If $X$ satisfies the hypotheses of Theorem {\rm \ref{realmain},}
  then $\{ K(n)_*(Q\Sigma ^rX) , r\geq 0\}$ forms a
$Z^+$\/{\rm -}\/indexed family of algebras equipped with a 
Dyer\/{\rm -}\/Lashof length-like filtration{\rm .}
\endproclaim

\demo{Proof}
Define an increasing  filtration F on  $K(n)_*(Q\Sigma ^mX)$  by
$$F_i(K(n)_* (Q\Sigma ^mX))=
\Im   ( K(n)_*(\mathbold{\cup}_{j\leq  i}E\Sigma _j\times _{ \Sigma _j}X^j)\rightarrow
    K(n)_*(Q\Sigma ^mX)).$$
The assertions (i), (iii), and (iv) will be proved  in Section 
\ref{kn;poly}.     Theorem
\ref{snsp} implies (vii),    (ii) follows from
 the property of the homology suspension, and (v), (vi), and (ix) are obvious from 
the definition.  Thus it remains to show (viii).  

This seems to be well-known to experts,
but does not seem to be in the published literature,  so we record a proof here.
According to Proposition 5.2 of 
\cite{geo},  one has the following commutative diagram where $C_nX$ denotes a
certain combinatorial model for $\Omega ^n\Sigma ^nX$. 
\figin{fig1}{1000}%
\noindent 
Furthermore,  the right vertical arrow factors as 
$$C_{n+1}X\stackrel{\beta _n}{\rightarrow }\Omega C_n\Sigma ^nX \stackrel{
\gamma _n}{\rightarrow } \Omega ^{n+1}\Sigma ^{n+1}X$$ according to 
[28, Prop.\ 5.4].  Its proof indicates that $\beta _n$'s and
$\gamma _n$'s are compatible with inclusions $C_nX\hra C_{n+1}X$, $\Omega
C_n\Sigma X\hra \Omega C_{n+1}\Sigma X$, and $\Omega ^n\Sigma ^nX \hra
\Omega ^{n+1}\Sigma ^{n+1}X$,  and that $\beta _n(F_lX)\subset
\Omega F_l(\Sigma X)$. Therefore one can pass to the colimit to get
 a commutative  
diagram
\figin{fig2b}{1000}%
\noindent 
such that $\beta (F_l(X))\subset 
\Omega (F_l(\Sigma X))$.  Now, forgetting the $CX$ on the upper-right corner and
 by taking the adjoint,  one gets the following
commutative diagram: 
\figin{fig3}{1000}%
\noindent 
with $j(\Sigma F_l(X))\subset F_l(\Sigma X)$.  Since the homology suspension 
map $h_*(QX)\rightarrow h_{*+1}(Q\Sigma X)$ is the composition:
$$h_*(QX)\rightarrow h_{*+1}(\Sigma QX) \stackrel{\rho _*}{\rightarrow }
h_*(Q\Sigma X),$$ one gets the desired result.\enddemo 

\proclaim{Proposition}
Let $X$ be a {\rm ($-1$)-}\/connected spectrum whose stable cells are all in even degrees{\rm .}
Then $K(n)_*(\und{X}{i})$ ($i\in Z^+$) forms a family of algebras with 
Dyer\/{\rm -}\/Lashof length-like filtration.
In particular{\rm ,}  this applies to the case when $X={\rm BP}${\rm .}
\endproclaim

\demo{Proof}
According to \cite{bpqs} there is an isomorphism of algebras
$K(n)_*(\und{X}{i})\cong \otimes K(n)^*(QS^j)$
 corresponding to a stable cellular
decomposition of $X$. Although there is nothing canonical in this decomposition,
once one fixes the decomposition for $K(n)_*(\und{X}{0})$,   one can choose the
rest to be compatible with the suspension homomorphism.  Thus we get the desired
result using  Propositions \ref{dll;ten} and \ref{dll;col}.\enddemo

Now we state a very useful property of  Dyer-Lashof length-like filtrations.

\proclaim{Proposition} \label{dll;main}
 Let $A$ be a $Z_{\geq d}$\/{\rm -}\/indexed family of algebras equipped with a
 Dyer\/{\rm -}\/Lashof length-like filtration{\rm .}  Let $B$ be another 
$Z_{\geq d}$\/{\rm -}\/indexed family of algebras equipped with suspension maps 
satisfying the properties {\rm (i)--(iii)} in   Definition {\rm \ref{filt;dll},}  such that
$\sigma $ sends the exterior part into the polynomial part and vice versa{\rm .}
Suppose that $f_i:A_i\ra B_i$ is a homomorphism such that $f_i$\/{\rm '}\/s commute
with the suspension and they respect the tensor product decomposition of
{\rm (iii).}  If 
$${\colim}_{i} f_i:\colim (A_{i,*}\stackrel{\sigma }{\ra }
A_{i+1,*+1}\stackrel{\sigma }{\ra }\cdots )
\ra
\colim (B_{i,*}\stackrel{\sigma }{\ra }
B_{i+1,*+1}\stackrel{\sigma }{\ra }\cdots )$$ is a monomorphism{\rm ,}  then so is
${\rm Ind}(f_i):{\rm Ind}(A_i)\ra {\rm Ind}(B_i)$ for each $i${\rm .}
\endproclaim

\demo{Proof}
We prove by induction on $l$  that $F_{p^l}({\rm Ind}(A_i))$ injects to 
${\rm Ind}(B_i)$.  Condition 
(ix) combined with the assumption of the proposition provides the first step. 
From conditions (ii),  (iii) and (iv),  we deduce that if
$$x\in ( \Im\, \sigma:{\rm Ind}(E_{i,*})\ra P_{i+1,*+1})\,\cap\, \Ker (P_{i+1,*+1}\ra
QP_{i+1,*+1})$$ then there exists $y$ such that $x=y^{p^r}$,  and $y\not \in
\Ker (P_{i+1,*+1}\ra QP_{i+1,*+1})$.  Now suppose that $a\neq 0$ is in
the kernel of the map $F_{p^l}({\rm Ind}(A_i))\ra {\rm Ind}(B)$.  As 
 ${\colim }_{i}f_i$
is a monomorphism,  there exists $r >0$ such that  $\sigma ^r(a)\neq 0$, $
\sigma ^{r+1}(a)=0$.  Thus the arguments above show that there exists $b$ such that
$\sigma ^r(a)=b^{p^s}$ and  that
$b$ reduces nontrivially to ${\rm Ind}(A_{i+r})$.  By (vii) we can conclude that $b\in F_{p^{l-s}}(A)$.
  However, by the  induction hypothesis,
$b$ maps nontrivially to ${\rm Ind}(B_i)$.  Since $b$ has to be primitive this means
that $f_i(b)$ is in the exterior part of $B_i$,  which means that $\sigma 
(f_i(b))$ is in the polynomial part.  However,  since $\sigma (b)$ is in the exterior
part,    $\sigma (f_i(b))$ has to be trivial. This contradicts 
condition
(iv).  \enddemo

\demo{Remark} This is essentially  how it was shown in \cite{W} that 
the subalgebra of $H_*(QS^i)$ mentioned in Remark
\ref{dll;name} injects to $H_*(\und{\rm BP}{i})$.
\enddemo

\section{How to recover ${\rm BP}$-cohomology from Morava $K$-theory}
\label{rwy;rev}

In this section we recall relevant results in \cite{rwy} and generalize them to
suit  our purpose.  There exist generalized cohomology 
theories $E(k,n)$ and $P(n)$ with the ring of coefficients 
$E(k,n)_*\cong E(n)_*/(p,v_1,\dots v_{n-1})$, $P(n)_*\cong {\rm BP}_*/(p,v_1,\dots
v_{n-1})$.  

\proclaimtitle{\cite{rwy}}
\proclaim{Lemma}
If $K(n)^{\rm odd}(X)=0$ then $E(k,n)^{\rm odd}(X)=0$ for $0\leq k\leq n${\rm ,}  and
$E(k,n)^*(X)$ has no $v_k$\/{\rm -}\/torsion{\rm .}
\endproclaim

\proclaimtitle{\cite{rwy}}
\proclaim{Theorem}
\label{rwy;ten}
Consider the following conditions{\rm .}
\medbreak
{\rm \listitem{i}} $K(n)^{\rm odd}(X)=0$ for an infinite number of $n$.
\medbreak{\rm \listitem{ii}} $E(k,n)^{\rm odd}(X)=0$ $\forall k,0<k<n$ for an infinite number of $n${\rm .}
\medbreak{\rm \listitem{iii}} $P(k)^{\rm odd}(X)=0$ for all $k${\rm .}
\medbreak{\rm \listitem{iv}} $K(k)^{\rm odd}(X)=0$ for all $k${\rm .}
\medbreak{\rm \listitem{v}} $E(k,n)^*(X)$ is $v_k$\/{\rm -}\/torsion free for all $0<k<n${\rm .}
\medbreak{\rm \listitem{vi}} $P(k)^*(X)$ is $v_k$\/{\rm -}\/torsion free for all $k${\rm .}
\medbreak{\rm \listitem{vii}} $(p,v_1,v_2,\dots ,)$ is a regular sequence in ${\rm BP}^*(X)${\rm .}
\medbreak{\rm \listitem{viii}} ${\rm BP}^*(X)\hat{\otimes }_{{\rm BP}^*}P(k)^*$ is isomorphic to $P(k)^*(X)$ for
all $k${\rm .}
\medbreak{\rm \listitem{ix}} ${\rm BP}^*(X)\hat{\otimes }_{{\rm BP}^*}K(k)^*$ is isomorphic to $K(k)^*(X)$ for
all $k${\rm .}
\medbreak{\rm \listitem{x}} ${\rm BP}^*(X)\hat{\otimes }_{{\rm BP}^*}K(k)^*$ surjects to $K(k)^*(X)$ for
all $k${\rm .}
\medbreak{\rm \listitem{xi}} ${\rm BP}^*(X)\hat{\otimes }_{{\rm BP}^*}E(k,n)^*$ surjects to $E(k,n)^*(X)$ for 
any $n\geq k\geq 0${\rm .}
\medbreak\noindent 
The conditions from {\rm (i)} to {\rm (iv)} are equivalent{\rm ,}  and they imply the rest{\rm .}
\endproclaim

\proclaimtitle{\cite{rwy}}
\proclaim{Theorem} \label{rwy;exa}
Let $X_i${\rm ,} $i=1,2,3$ be spaces satisfying one of the conditions from {\rm (i)} to {\rm (iv)}
of the theorem above{\rm .}  Let $f:X_1\rightarrow X_2${\rm ,} $g:X_2\rightarrow X_3$ be maps
with $g\circ f$ null\/{\rm -}\/homotopic{\rm .}  Then{\rm ,}
\medbreak 
{\rm \listitem{i}} If $K(n)^*(g)$ is mono for all $n$ then so is ${\rm BP}^*(g)${\rm .}
\medbreak {\rm \listitem{ii}} If $K(n)^*(f)$ is epi  for all $n$ then so is ${\rm BP}^*(f)${\rm .}
\medbreak {\rm \listitem{iii}} Furthermore if all the spaces are $H$\/{\rm -}\/spaces and maps are $H$\/{\rm -}\/maps
such that 
$$K(n)_*(X_1)\stackrel{K(n)_*(f)}{\rightarrow }K(n)_*(X_2)
\stackrel{K(n)_*(g)}{\rightarrow }K(n)_*(X_3)$$
is an exact sequence of Hopf algebras{\rm ,}  then 
$${\rm BP}^*(X_3)\stackrel{{\rm BP}^*(g)}{\rightarrow }{\rm BP}^*(X_2)
\stackrel{{\rm BP}^*(f)}{\rightarrow }{\rm BP}_*(X_1)$$ is a coexact sequence of augmented
${\rm BP}^*$-algebras{\rm .}  That is{\rm ,} ${\rm BP}^*(f)$ is the cokernel of the map ${\rm BP}^*(g)$ in
the category of augmented ${\rm BP}^*$\/{\rm -}\/algebras{\rm .}  More concretely{\rm ,}  ${\rm BP}^*(g)$ induces
an isomorphism between the quotient of ${\rm BP}^*(X_2)$ by the image of the 
augmentation ideal of ${\rm BP}^*(X_3)$ by ${\rm BP}^*(f)$ and ${\rm BP}^*(X_1)${\rm .}
\endproclaim

\proclaimtitle{\cite{rwy}}
\proclaim{Theorem} \label{rwy;lan}
Let $X${\rm ,} $Y$ be spaces satisfying one of the conditions from {\rm (i)} to {\rm (iv).}  Then
$${\rm BP}^*(X\times Y)\cong {\rm BP}^*(X)\hat{\otimes }_{{\rm BP}^*}{\rm BP}^*(Y).$$
\endproclaim

\demo{Remark} Note that this does not follow from Theorem \ref{rwy;ten} and
Landweber's exact functor theorem unless $Y$ is finite.  A naive ``proof" 
would involve commuting a direct limit with an inverse limit.\enddemo

Now we generalize these results.

\proclaim{Theorem} \label{land;gen}
The conditions from {\rm (vi)} to {\rm (x)} in   Theorem {\rm \ref{rwy;ten}} are equivalent{\rm .}
They are also equivalent to {\rm (v)} with $0<k<n$ replaced by $0\leq k<n$ {\rm (}\/called condition {\rm (v)$'$).}
Furthermore{\rm ,}  it suffices to assume one of these equivalent conditions on 
spaces appearing in  Theorems {\rm \ref{rwy;exa}} and {\rm \ref{rwy;lan}} to obtain the
same conclusion{\rm .}
\endproclaim

\demo{Proof} It is well-known that   condition (vi) implies (vii) that implies (viii)
(this can be shown easily, inductively, from the cofibration sequence 
$P(k)\stackrel{v_k}{\rightarrow } P(k) \rightarrow P(k+1))$.
In \cite{rwy} it was shown that $P(k)^*(X)\hra \mathbold{\oplus}_{n>k} E(k,n)^*(X)$,
so that (v)$'$ implies (vi).  Since by   Morava's little structure theorem (\cite{jw},
\cite[Prop.\ 1.9]{w2} for the present form)
$$P(k)^*(X)\hat{\otimes }_{P(k)^*}(K(k)^*)\cong
K(k)^*(X),$$ (viii) implies (ix) which obviously implies (x).  The cofibration
sequence $$E(k,n)\stackrel{v_k}{\rightarrow } E(k,n)\rightarrow E(k+1,n)$$
and the fact that the map from ${\rm BP}$ to
$E(k+1,n)$ factors through $E(k,n)$ can be used to show that (xi) implies (v)$'$.
The same cofibration sequence and the fact that the filtration by the power of
$v_k$ is complete in ${\rm BP}^*(X)$ prove that (x) implies (xi) by downward induction
on $k$,  where (x) serves as the starting point of the induction.  This finishes
the proof of the equivalence of the conditions listed.  The 
proof of  Theorems \ref{rwy;exa} and \ref{rwy;lan} does not really rely on the
properties from (i) to (iv),  but  only uses the properties (v) and (vi) (and 
other properties that hold for ${\rm BP}^*(X)$ for any space $X$);  more precisely
the fact that the long exact sequences associated to the aforementioned 
cofibrations become just a bunch of short exact sequences.  Therefore it suffices
to assume one of these conditions to get the same conclusion.  \enddemo

\section{${\rm BP}$-cohomology of the extended power construction} \label{wreath}

Morava $K$-theory of the extended power construction was first studied in
\cite{Hu1}, \cite{HKR}.  The work in \cite{jrh;gen} treats the most general  
situation,  as well as it deals with the case of other complex oriented 
cohomology theories including ${\rm BP}$-cohomology.  We use the results in 
\cite{jrh;gen} to obtain:

\proclaim{Theorem} \label{th51}
Let $X$ be a connected space satisfying one of the equivalent conditions in Theorem 
{\rm \ref{land;gen}.}  Then $D_{Z/p}X$ satisfies the same conditions.{\rm }  Furthermore{\rm ,}
if ${\rm BP}^*(X)$ is well\/{\rm -}\/generated then so is ${\rm BP}^*(D_{Z/p}X)${\rm .}
  Here $D_{Z/p}X=EZ/p\times _{Z/p}X^p$
where 
$Z/p$ acts on $X^p$ by permutation{\rm .}
\endproclaim

Before proving the theorem,  we recall a result on the behavior of the
Atiyah-Hirzebruch-Serre spectral sequence for the fibration
$X^p\ra D_{Z/p}X\ra BZ/p$. 
First note that,  if $h$ is a generalized cohomology theory,  the 
Atiyah-Hirzebruch spectral sequence for the space $BZ/p$ acts on the
AHSss in question.  
We only consider the case when 
$h^*( X^p)$ is isomorphic to $ h^*(X)^{\cotimes p}$.
We say that the AHSss
$H^*(BZ/p,h^*(X^p))\Rightarrow h^*(D_{Z/p}X)$ is simple,  if there is no
differential in this spectral sequence other than those that are forced by 
the action of the AHSss $H^*(BZ/p,h^*)\Rightarrow h^*(BZ/p)$. 

More precisely,

\demo{Definition {\rm 4.2 (\cite{jrh;gen})}}
If $h^*( X^p)\cong h^*(X)^{\cotimes p}$, we say that the AHSss 
$E_2^{*,*}\cong H^*(BZ/p,h^*(X^p))\Rightarrow E_{\infty }^{*,*}\cong 
h^*(D_{Z/p}X)$ is {\it simple} if $E_2^{0,*}\cong E_{\infty }^{0,*}$.
\enddemo
\advance\theoremcount by 1

We need to know the behavior of this AHSss in more detail.  
The $E_2$ term is isomorphic to $A^*\otimes H^*(BZ/p,h^*)\oplus B^*$,  where
$A_*$ is the span of the elements of the form $\stackrel{p\mbox{ factors}}{
\overbrace{a \otimes \cdots \otimes a}} $,  $a\in h^*(X)$,  whereas 
$B^*$ is the span of the elements of the form $\Sigma _{\sigma \in Z/p}
\sigma (a_1\otimes \cdots \otimes a_p)$.

\demo{Proof of Theorem {\rm \ref{th51}}}
According to Theorems 2.5 (or the remark preceeding  it) and 6.1 of
\cite{jrh;gen},  the condition for $X$ implies that both 
the AHSss for ${\rm BP}^*(D_{Z/p}X)$ and $K(n)^*(D_{Z/p}X)$ are simple.  
This means for $K(n)$,  $E_{\infty }^{*,*}\cong A^*\otimes H_*(H^*(BZ/p,K(n)_*),
v_nQ_n)\oplus B^*$,  where $A^*$ and $B^*$ are as above,  and $Q_n$ is the $n^{\rm th}$
Milnor's Bockstein operation.   Thus as an algebra over $K(n)_*$,  it is generated by the
elements of
$A^*$,  $B^*$,  and the element in $E_{\infty }^{2,*}$ represented by 
the element $0\neq x\in H^2(BZ/p,K(n)_*)$.  The collapsing of its ${\rm BP}$-counterpart
$E{'}_*^{*,*}$ (which is nothing but the simpleness for ${\rm BP}$)
implies that all these elements are in the image of the map
$E{'}_{\infty }^{*,*}\ra E_{\infty }^{*,*}$ induced by the natural transformation
${\rm BP}^*(-)\ra K(n)^*(-)$ up to multiplication by some power of $v_n$.
Thus the 
condition (x) of Theorem
\ref{land;gen} is easily seen to be satisfied for $D_{Z/p}X$.  The second 
statement follows immediately from the collapsing of the 
 AHSss for ${\rm BP}^*(D_{Z/p}X)$.
\enddemo

As usual,  properties that are preserved by the construction $D_{Z/p}$ are
preserved by the construction $Q$.  Namely;

\demo{Proof of Theorem {\rm \ref{qx;ten}}}
Since these two properties only concern   the\break ${\rm BP}^*$-module structure,  by
Theorem \ref{snsp} it suffices to show these properties for $E\Sigma _n^+ \wedge
X^n$.  However,  one can easily show by transfer arguments that\break $p$-locally,
these spaces are stable retracts of products of the spaces of the form 
$D_{Z/p}(\cdots (D_{Z/p}(X)))$ (see, e.g. \cite{mc}). 
 Thus one obtains the desired result from
  Theorem \ref{th51}. \enddemo

\section{The image of the Thom map} \label{thom;im}

In this section we describe $M_{QX}$ in terms of $M_X$
with some hypotheses on $X$.   First we establish an upper bound on $M_{QX}$.

\proclaim{Proposition} \label{bpg;det}
Let $X$ be a connected space with the property ${\rm BP}^*(X^j)\break \cong {\rm BP}^*(X)^{\hat{\otimes }j}${\rm .}
Let $$B=\{f|f\in H_*(X)=\Hom (H^*(X);Z/p), \mbox{$f$ vanishes on $M_X$}.\}$$  Choose
its complement $A${\rm ,}  i.e.{\rm ,} a subspace of $H_*(X)$ such that
$H_*(X)=A\oplus B${\rm .}  Then one has
$$M_{QX}\subset C=\{ \phi |H_*(QX)\rightarrow Z/p
\mbox{ such that } \phi (S)=0\} $$
where $S$ is the ideal generated by the elements of the form 
$Q^Ix$ with $x\in B$ or of the form $Q^Jx$ with $x\in A$ and $J$ containing at
least one Bockstein{\rm ,}  and $A$ and $B$ are considered as subspaces of $H_*(QX)$ via
the inclusion $H_*(X)\hra H_*(QX)${\rm .}
\endproclaim

\demo{Proof} This will be proved in three steps.  First we prove that
the elements of the form $Q^Ix$ with $x\in B$ or of the form 
$Q^Ix$ with $x\in A$ and $I$ containing at
least one Bockstein can be written 
as  a linear combination of elements of the form
$Q^{ \circ K}[1]\circ z$ with 
either $K$ containing at
least one Bockstein or $z\in B$,
where $Q^{\circ  K}[1]$ denotes the 
element $\beta ^{\varepsilon _1}Q^{s_1}[1]\circ \cdots \circ \beta ^{\varepsilon _l}
Q^{s_l}[1]$ if $K=\varepsilon _1,s_1,\dots \varepsilon _l,s_l$. 
  We prove the  following two statements by induction on 
$l(I)$ and ${\rm deg}(Q^{I'}(x))$, where if $I=(\varepsilon ^1,I^1,\varepsilon ^2,I^2,\dots)$ then $I'=(\varepsilon
^2,I^2,\dots)$.
\medbreak
\listitem{i} Let $N\subset H_*(X)$ be a subspace closed under the action of the Steenrod algebra.
Then $Q^Ix$ ($x\in N$) can be written as a sum of the elements of the form
$Q^{\circ K}[1]\circ z$ with $z\in N$.
\medbreak
\listitem{ii} Furthermore suppose that $\beta ( H_*(X))\subset N$.
If $I$ contains a Bockstein,  then $Q^Ix$ ($x\in H_*(X)$) can be
 written as a sum of the elements of the form $Q^{\circ K}[1]\circ z$ with 
either $K$ containing a Bockstein or $z\in N$.
\medbreak
To prove the first statement, using May's formula,  one gets
$$Q^I(x)=\beta ^{\varepsilon ^1}Q^{I^1}[1]\circ Q^{I'}(x)-\Sigma _{t>0}\beta 
^{\varepsilon ^1}Q^{I^1+t}P^t_*Q^{I'}(x)+\varepsilon ^1\Sigma _{t\geq 0}Q^{I^1+t}P^t_*
\beta Q^{I'}(x).$$  Since $l(I')=l(I)-1$,  the first term can be taken care of
by  induction on $l$.  Using Nishida relations one can rewrite 
$P^t_*Q^{I'}(x)$ ($t>0$) and $P^t_*\beta Q^{I'}(x)$ as a linear combination of
elements of the form $Q^J(z)$ with $l(J)=l(I')$,  $z$ belonging to the orbit of
$x$ by the action of the Steenrod algebra,  thus belonging  to $N$ and ${\rm deg}(Q^J(z))<
{\rm deg}(Q^{I'}(x))$.  Since the sequences
$(\varepsilon ^1,I_1,J)$ have the same length as $I$,  and the degree of $Q^J(z)$
is less than that of $Q^{I'}(x)$,  the two summations can be taken care of by
the
induction hypothesis,  which finishes the proof of~(i).  To prove (ii),  when
$\varepsilon ^1=0$,  $I'$ still contains a Bockstein,  and the rest of the argument
is similar.  When $\varepsilon _1=1$,  one can treat the terms in the two 
summations similarly,  and,  to take care of
the term $\beta ^{\varepsilon _1}Q^{I^1}[1]\circ Q^{I'}(x)$,  one applies the case
$N=H_*(X)$ of (i) to $Q^{I'}(x)$.  Thus one has proved (i)\break and~(ii).

Now note that since the operations $P^j$s are covered by some Landweber-Novikov
operations,  $B$ is stable under the action of the  $P^j_*$'s,  and since the
Bockstein vanishes on $M_X$,  in $H_*(X)$ the image of $\beta $ is contained in
$B$.  Thus one can take $N$ in the statements above to be $B$ to obtain the
desired result.

In the next step,  we show that
there exists a family of spaces $Y_{k,l,i}\cong (B\Sigma _p)^l
\times X^k$
and a map $g_{k,l,i}:Y_{k,l,i} \rightarrow QX$ such that any element of $S$ can
be written as a linear combination of elements of the form
 $g_{k,l,i\mbox{ } *}(e_{j_1}\otimes \cdots \otimes e_{j_l}
\otimes x_1\otimes \cdots \otimes x_k)\in H_*(QX)
$ where either at least one of $j_t$'s is congruent to $(-1)$ mod $2(p-1)$ or 
at least one of the $x_k$'s is in $B$.

As a matter of fact it is enough to take the family 
\bea
% \lefteqn
& &\hskip-.75in {\{ \overbrace{(B\Sigma _p \times \cdots \times B\Sigma 
_p)}^{m_1\mbox{ factors}}\times X \times\cdots \times \overbrace{(B\Sigma _p \times \cdots
 \times B\Sigma _p)}^{m_k\mbox{ factors}} \times X}\\ & \rightarrow & B\Sigma _{p^{m_1}}\times X \times
\cdots \times B\Sigma _{p^{m_k}}\times X\\ & \rightarrow &  
B\Sigma _{p^{m_1}+\cdots +p^{m_{k}}}\times X^k\\
 & \rightarrow & QS^0\times QX\\
 & \rightarrow & QX\} 
 \eea
where the first map is induced by the multiplication map $B\Sigma _a\times B\Sigma
 _b\rightarrow B\Sigma _{ab}$, the second by the addition $B\Sigma _a \times 
B\Sigma _b\rightarrow B\Sigma _{a+b}$,  the third by obvious ones,
and the last by the composition pairing.   Now by definitions 
$\alpha _1\otimes \cdots \otimes \alpha _k$   each
$\alpha _s\in H_*(B\Sigma _p 
^{m_s}\times X)$ is mapped to 
$\mu _*(\alpha _1)\star \cdots \star \mu _*(\alpha _k)$ where 
$\mu _*$ is given by 
$$\mu _*((e_{2j_1(p-1)+\varepsilon _1}\otimes \cdots \otimes e_{2j_{m}(p-1)+
\varepsilon _m})\otimes x )
=(\beta ^{\varepsilon _1}Q^{j_1}[1]\circ \cdots \circ \beta ^{\varepsilon _m} 
Q^{j_{m}}[1])\circ x
.$$  Thus one deduces the desired result from the previous step.

Now the proof of the proposition  can be completed as follows.
The conclusion is  equivalent to the vanishing of the restriction to
$M_{QX}\otimes S$ of  the Kronecker pairing $H^*(X)\otimes H_*(X)\ra Z/p$.
However, the assumption on $X$ implies that 
${\rm BP}^*(Y_{k,l,i})\cong {\rm BP}^*(B\Sigma _p)^{\cotimes l}\bpg {\rm BP}^*(X)^{\cotimes k}$.
Since we know that if $c\in M_{B\Sigma _p}$ then
$c(e_j)=0$ if $j$ is congruent to $(-1)$ mod $2(p-1)$,  we see that
if $c\in M_{Y_{k,l,i}}$, $ c$ vanishes on elements of the
form $e_{j_1}\otimes \cdots \otimes e_{j_l}
\otimes x_1\otimes \cdots \otimes x_k
$ where either at least one of $j_t$'s is congruent to $(-1)$ mod $2(p-1)$ or
at least one of the $x_k$'s is in $B$.  Thus if $f\in M_{QX}$, then  $g_{k,l,i}^*(f)$
vanishes on such elements.  The result of the previous step  now implies the
desired result.
\enddemo

Next we go on to establish a lower bound for $M_{QX}$.
\proclaim{Lemma} \label{tho;lb}
Let $\{ f_i:X\ra \und{\rm BP}{d_i}|i\in I\} $ be a set of topological
${\rm BP}^*{\rm (BP)}$\/{\rm -}\/module 
generators for $\wt{\rm BP}^*(X)${\rm .}  Then one has
$$\Im (\otimes _{i\in I}H^*(\Omega ^{\infty }f_i):
H^*(\Pi _{i\in I}\und{\rm BP}{d_i})\ra H^*(QX))=C$$
where $C$ is as in   Proposition {\rm \ref{bpg;det}.}
\endproclaim

\demo{Proof} Consider a $Z^+$-indexed family of graded algebras 
$\{ H_*(Q\Sigma ^nX)|n\break \in Z^+\}$.  We fix a direct sum decomposition as in
Proposition \ref{bpg;det}: $H_*(X)=A_X\oplus B_X$ and $H_*(\Sigma ^nX)
=A_{\Sigma ^nX}\oplus B_{\Sigma ^nX}$ compatible with the suspension isomorphism.
Let $T _{\Sigma ^nX}$ denote the subalgebra of $\ H_*(Q\Sigma ^nX)$ generated 
by the elements of the form $Q_Jx$ with $x\in A_{\Sigma ^nX}$ and $J$ containing
no Bockstein.  Then $\ H_*(Q\Sigma ^nX)=T _{\Sigma ^nX}\oplus S_{\Sigma ^nX}$,
where $S$ is as defined in Proposition \ref{bpg;det}.  It is easy to see that
$\{ T _{\Sigma ^nX}|n\in Z^+\}$ forms a $Z^+$-indexed family of graded algebras
with Dyer-Lashof length-like filtration.  Furthermore,  one has 
$$\colim (\ra T _{\Sigma ^nX}\ra T _{\Sigma ^{n+1}X}\ra \cdots )\cong A.
$$
On the other hand, by   choice of the maps $f_i$,  we know that
${\rm BP}^*(\vee _{i\in I}\Sigma ^{d_i}{\rm BP})$ surjects to $\wt{\rm BP}^*(X)$.   Thus
$H^*(\vee _{i\in I}\Sigma ^{d_i}{\rm BP})$ surjects to $M_X$ (note that the Thom 
homomorphism ${\rm BP}^*(\vee _{i\in I}\Sigma ^{d_i}{\rm BP})\ra H^*(\vee _{i\in I}\Sigma ^{d_i}{\rm BP})$ is also surjective so that $H^*(\vee _{i\in I}\Sigma ^{d_i}{\rm BP})$ maps
to $M_X$).  As
 the restriction of the pairing $H^*(X)$ and $H_*(X)$ identifies $A$
with the dual of $M$,  $A$ injects to $H_*(\vee _{i\in I}\Sigma ^{d_i}{\rm BP})$.
Thus by Theorem \ref{dll;main} we see that $T _X$ injects to 
$H_*(\Pi _{i\in I}\und{\rm BP}{d_i})$.  On the other hand, since 
${\rm BP}^*(\und{\rm BP}{d_i})$  surjects to $H^*(\und{\rm BP}{d_i})$,  $B_X$
is seen to be in the kernel of $H_*(\Omega ^{\infty }f_i)$.  Furthermore,  since Bocksteins act
trivially on $H_*(\und{\rm BP}{d_i})$ we see that $S_X\subset (\Ker \oplus _i
H_*(
\Omega ^{\infty }f_i))$.  As one has seen $H_*(QX)=S_X\oplus T_X$,   
this shows that
 $S_X=\Ker (\oplus _iH_*( \Omega ^{\infty }f_i))$.  By dualizing one gets 
the desired result.  \enddemo

\demo{Proof of Theorem {\rm \ref{int:qbpg}}} Since ${\rm BP}^*(\und{\rm BP}{d_i})$ surjects to $H^*(\und{\rm BP}{d_i})$,  
Lemma \ref{tho;lb} implies that $C\subset M_{QX}$.  Combining this with
Proposition \ref{bpg;det},  one gets the desired result. \enddemo

\section{$K(n)_*(QX)$ and ${\rm BP}^*(QX)$} \label{kn;poly}

In this section we determine $K(n)_*(QX)$ and ${\rm BP}^*(QX)$ for the spaces $X$
satisfying the hypotheses of Theorem \ref{realmain}.
First we prove:

\proclaim{Theorem} \label{kn;mon}
Let $X$ be a connected space satisfying one of the equivalent conditions of Theorem
{\rm \ref{land;gen},}
and $${\rm BP}^*(X)\bpg Z/p\hra H^*(X;Z/p).$$
Let $\{ f_i:X\ra \und{\rm BP}{d_i}|i\in I\} $ be a set of ${\rm BP}^*{\rm (BP)}$\/{\rm -}\/module 
generators for $\wt{\rm BP}^*(X)$\/{\rm .}\/  Then 
\medbreak
{\rm \listitem{i}} $\hat{\otimes } _i{\rm BP}^*(f_i)$ is surjective{\rm ,}
\medbreak
{\rm \listitem{ii}} $\otimes _iK(n)_*(f_i)$ is injective{\rm .}
\endproclaim

\demo{Proof}
According to Theorem \ref{land;gen},  one has  
${\rm BP}^*(X^j)\cong {\rm BP}^*(X)^{\cotimes j}$.  Thus we see that $M_{QX}$ agrees with the
image of the composition $${\rm BP}^*(\Pi _{i\in I}\und{\rm BP}{d_i}) \ra {\rm BP}^*(QX)\ra H^*(QX).
$$
  However,  we also
know from Theorem \ref{qx;ten} that $${\rm BP}^*(QX)\bpg (Z/p)\subset H^*(QX).$$  This
concludes the proof of (i). Using Theorem \ref{qx;ten} one can deduce (ii) from~(i). 
\enddemo

Note that in the above,  one can take all $d_i$ to be positive.  This easily 
follows from the fact that ${\rm BP}^*(X)$ is well-generated.  One can also deduce it
by the theorem of Quillen \cite{Q} which says that the ${\rm BP}$ cohomology of a space
is always generated by nonnegative degree elements.  Thus
 $K(n)_*(\Pi _{i\in I}\und{\rm BP}{d_i})$
 is a free commutative algebra, and by Theorem B.7 of
\cite{BKf} any of its Hopf subalgebras is a free commutative algebra.
Thus we obtain the first statement of  
Theorem \ref{int:knfree}.  Now we will go on to show the second statement.
\demo{Proof of Theorem \ref{int:knfree}}
We will study the bar spectral sequence for the fibration 
$QX\ra {\rm pt} \ra Q\Sigma X$. 
We have just seen that there is a Hopf algebra isomorphism $K(n)_*QX\cong P_X
\otimes E_X$, where $P_X$ is a polynomial algebra concentrated in even degrees
and $E_X$ is an exterior algebra generated by odd degree elements.
 Thus we have 
\bea
E^2 & \cong &  {\rm Tor}^{K(n)_*(QX)}(K(n)_*,K(n)_*)
\\ & \cong & \Gamma (\sigma {\rm Ind}(E_X)) \otimes \Lambda (\sigma {\rm Ind}(P_X)).
\eea
Now,  let $J=\{i\in I|d_i \mbox{ is even}\}$,  and let $Y$ be the cofiber
of the map $f'=\Pi _{i\in J}f_i:X\ra \vee _{i\in J}\Sigma ^{d_i}{\rm BP}$.  
Consider the bar spectral sequence associated to the fibration 
$QX \ra \Pi _{i\in J} \und{\rm BP}{d_i} \ra \Omega ^{\infty }Y$.  As the map
$K(n)_*(\Pi _{i\in I} \Omega (f_i))$ is injective,  we see that the map
$K(n)_*(\Pi _{i\in J} \Omega (f_i))$ is injective on the
polynomial part.  Thus
the $E^2$ term has  the form $$ {\rm Tor}^{E_X}(K(n)_*,K(n)_*) \otimes
K(n)_*(\Pi _{i\in J} \und{\rm BP}{d_i})//P_X,$$
 which is concentrated in even 
degrees.  Therefore it collapses.  Next we compare it with 
the bar spectral sequence
above,  and we see that the factor $\Gamma (\sigma {\rm Ind}(E_X))$  is composed of
permanent cycles only.  As the other factor $\Lambda (\sigma {\rm Ind}(P_X))$ is also
composed of permanent cycles since it is generated by homological degree 1 
elements,  we see that this spectral sequence collapses as well.  Thus $E^2=
E^{\infty }$ and the $E^{\infty }$ term is a cofree coalgebra; thus there can be
no coalgebra extension and $K(n)_*(Q\Sigma X)$ is a cofree coalgebra.
\enddemo

Note that   Proposition 5.1 above immediately implies the properties i) and iii) 
in Definition 2.2.  The remaining property iv) follows easily from
the collapse of the bar spectral sequence for the fibration
$QX \ra {\rm pt} \ra Q\Sigma X$ which  has just been proved.  This completes the proof of
Proposition \ref{dll;qx}.  We also note some special cases of   Proposition 5.1
above.

\proclaim{Corollary} \label{kn;biex} \hskip-8pt
Let $X$ be a connected space with $\widetilde{K(n)}_{\rm even}(X)=0$ for all~$n${\rm .}  Then 
$K(n)_*(QX)$ is an exterior Hopf algebra{\rm .}
\endproclaim

{\it Proof}. Since ${K(n)}_{\rm odd}(\Sigma X)=0$ for all $n$,   by Theorem 
\ref{rwy;ten}  one sees
that  ${\rm BP}^{{\rm odd}}\raise6pt\hbox{$\phantom{\int}$}\hskip-3pt (\Sigma X)=0$; i.e., $\widetilde{\rm BP}^{\rm
even}(X)=0$. Thus there is an inclusion of 
Hopf algebras $K(n)_*(QX)\ra K(n)_*(\Pi _i\und{\rm BP}{d_i})$ where all $d_i$'s are
odd.  Since $K(n)_*(\und{\rm BP}{\rm odd})$ is an exterior Hopf algebra,  we get the desired 
result. \hfill\qed

\proclaim{Corollary} Let $X$ be a connected space with $K(n)_{\rm odd}(X)=0$ for all~$n${\rm .}  Then
 $K(n)_*(QX)$ is a polynomial algebra{\rm .}  \endproclaim

{\it Proof}. The same arguments as above work  except that now all $d_i$'s are even.
\hfill\qed\smallbreak

For a  fixed value of $n$,   fewer $f_i$'s will suffice in Theorem \ref{kn;mon},  namely:

\proclaim{Proposition} \label{kn;mono}
Let 
 $\{ f_i:X\ra \und{\rm BP}{d_i}|i\in I\} $ be 
maps with the property that
$\widetilde{K(n)_*}(X)\stackrel{\oplus K(n)_*(f_i)}{\lrar}
\oplus \Sigma ^{d_i}K(n)_*{\rm (BP)}$ is injective{\rm .}
  Then the map
$${\rm Ind}(\otimes _iK(n)_*(f_i)):{\rm Ind}(K(n)_*(QX))\ra {\rm Ind}(K(n)_*(\Pi _i\und{\rm BP}{d_i}))
$$ is injective.
\endproclaim

{\it Proof}.
Now that one has seen that $K(n)_*(Q\Sigma ^rX)$'s are free,  the Proposition
\ref{dll;qx} implies that it is equipped with the Dyer-Lashof length-like 
filtration.  Thus we get the result by applying   Proposition \ref{dll;main}.
\hfill\qed\medbreak

We note another variant which should be of independent interest, namely:
\proclaim{Corollary}
Let $\{ f_i:X\ra \und{E(n)}{d_i}|i\in I\} $ be
maps with the property that
$\widetilde{K(n)_*}(X)\stackrel{\oplus K(n)_*(f_i)}{\lrar}
\oplus \Sigma ^{d_i}K(n)_*(E(n))$ is injective{\rm .}
  Then
the map
$${\rm Ind}(\otimes _iK(n)_*(f_i)):{\rm Ind}(K(n)_*(QX))\ra {\rm Ind}(K(n)_*(\Pi _i\und{E(n)}{d_i}))
$$ is injective{\rm .}
\endproclaim

 \demo{Proof} This follows from the fact that $K(n)_*(\und{E(n)}{i})$ is a
polynomial algebra if $i$ is even and an exterior algebra if $i$ is odd 
\cite{huad}, \cite{hh}.$\phantom{\sum^\int}$ \enddemo

Notably,  if we take $X$ to be a sphere,  we can use the unit map for the 
spectrum $E(n)$.  This generalizes the well-known result on 
injections $K(1)_*(QS^0)\hra K(1)_*(BU\times Z)$ and
$K(1)_*(QS^2)\hra K(1)_*(BU)$ (\cite{hs}).  (Strictly speaking,  the
case for $QS^0$ is not covered by the fact that  $S^0$ is  not   connected,  though
it is   not difficult to extend our result to this case,  which is left as an
exercise for interested readers.)

Regard $f_i$'s as maps of spectra $\Sigma ^{\infty }X\ra \Sigma ^{d_i}{\rm BP}$.
Let $C_f$ denote the cofiber of the map $f=\vee f_i:X\ra 
\vee \Sigma ^{d_i}{\rm BP}$.  We will now consider $K(n)_*(\Omega ^{\infty }C_f)$.

\proclaim{Lemma} \label{res;debut}
There is a
short exact sequence of Hopf algebras $$K(n)_*(QX)\ra K(n)_*(\Pi _i
\und{\rm BP}{d_i}) \ra K(n)_*(\Omega ^{\infty }C_f).
$$
\endproclaim

{\it Proof}. We consider the bar spectral sequence associated to the fibration 
$QX\ra \Pi _i \und{\rm BP}{d_i} \ra \Omega ^{\infty } C_f$.  
By  Proposition \ref{kn;mono} and Theorem 10.8 of \cite{BK},
we see that the $E_2$ term $ {\rm Tor}^{K(n)_*(QX)}( K(n)_*,K(n)_*(\Pi _i
\und{\rm BP}{d_i}))$ is concentrated in homological degree zero and isomorphic to
$K(n)_*(\Pi _i \und{\rm BP}{d_i})\otimes _{K(n)_*(QX)} K(n)_*$
 so that the SS
collapses and we get the desired result.\hfill\qed

\medbreak {\it Remark} 6.4. Although we can prove the injection 
${\rm Ind}(K(n)_*(QX)) \hra {\rm Ind}(K(n)_*(\Pi _i
\und{\rm BP}{d_i}))$,
%the remark made after the Lemma \ref{hop:spl} implies that the informations 
%obtained so far are not sufficient to say more about
% $K(n)_*(\Omega ^{\infty }C_f)$. 
as we are dealing with Hopf algebras with periodic gradings,  it does not
suffice to conclude that $K(n)_*(\Omega ^{\infty }C_f)$ is free. \advance\theoremcount by 1

\proclaim{Proposition} \label{cof;pol}
$K(n)_*(\Omega ^{\infty } C_f)$ is a free commutative algebra{\rm .}
\endproclaim

{\it Proof}. We use the notation  and definitions   in the proof of
Theorem \ref{int:knfree}.  The short exact sequence above splits as the 
tensor product of the short exact sequences 
$E_X\ra K(n)_*(\Pi _{i\in I-J}\und{\rm BP}{d_i}) \ra E_C$, $P_X\ra 
 K(n)_*(\Pi _{i\in J}\und{\rm BP}{d_i})
\ra P_C$ with $E_C\otimes P_C\cong K(n)_*(\Omega ^{\infty } C_f)$.  Obviously, $E_C$ is
 an exterior algebra generated by odd degree elements, and $P_C$
is concentrated in even degrees.  Thus it suffices to show that $P_C$ is a
polynomial algebra.  However, we have seen in the  proof of
Theorem \ref{int:knfree} that the bar spectral sequence for the fibration
$QX\ra \Pi _{i\in J}\und{\rm BP}{d_i} \ra \Omega ^{\infty }Y$ collapses,
which implies that $P_C$ injects to $K(n)_*(\Omega ^{\infty }Y)$.  Thus it 
suffices to prove that $K(n)_*(\Omega ^{\infty }Y)$ is a polynomial algebra.\

Consider the Eilenberg-Moore spectral sequences (see \cite{Lsm})
for the following three fibrations: $QX \ra {\rm pt} \ra Q\Sigma X$, $\Pi _{i\in J}\und{\rm BP}{d_i} \ra {\rm pt} \ra
\Pi _{i\in J}\und{\rm BP}{d_i+1}$, and $\Omega ^{\infty }Y\ra Q\Sigma X \ra
\Pi _{i\in J}\und{\rm BP}{d_i+1}$.  We show that they collapse at $E_2$ and   
actually converge.  We will call the
$E_2$ (thus $E_{\infty }$) term   $E_2(1), E_2(2), E_2(3)$\break  
respectively.   For the first fibration,  we have seen that
$K(n)_*(Q\Sigma X)\cong P_{\Sigma X}\break\otimes E_{\Sigma X}$, where 
$P_{\Sigma X}$ is isomorphic to $\Gamma (\sigma ({\rm Ind}(E_X)))$ as coalgebras, and\break
$ E_{\Sigma X} \cong \Lambda (\sigma ({\rm Ind}(P_X)))$.  Thus the $E_2$ term is 
\bea
%\noalign{\vskip-18pt}
 E_2(1) & \cong & 
{\rm Cotor} _{\Gamma (\sigma ({\rm Ind}(E_X))) \otimes \Lambda (\sigma ({\rm Ind}(P_X)))}(
K(n)_*,K(n)_*) \\ \noalign{\vskip4pt}
& \cong & E_X\otimes P_X \\ \noalign{\vskip4pt}
&\cong & K(n)_*(QX) \eea 
thus it collapses at $E_2$ and   converges. (Tamaki also constructed an
Eilenberg-Moore type spectral sequence that is strongly convergent
\cite{Ta}.)   Similarly the second collapses and converges.  As to the third one,
we see  that 
\bea E_2(3) & \cong &  {\rm Cotor} _{K(n)_*(\Pi _{i\in J}\und{\rm BP}{d_i+1}})(K(n)_*(Q\Sigma X),K(n)_*)
 \\ \noalign{\vskip4pt} 
& \cong & {\rm Cotor} _{K(n)_*(\Pi _{i\in J}\und{\rm BP}{d_i+1}}) (P_{\Sigma X} \otimes
E_{\Sigma X},K(n)_*) \\ \noalign{\vskip4pt}
& \cong & P_{\Sigma X}\otimes {\rm Cotor} _{K(n)_*(\Pi _{i\in J}\und{\rm BP}{d_i+1})//
E_{\Sigma X}}(K(n)_*,K(n)_*) \\ \noalign{\vskip4pt}
& \cong & P_{\Sigma X}\otimes {\rm Sym} (\sigma ^{-1} {\rm Ind} (K(n)_*(\Pi _{i\in J}\und{\rm BP}{d_i+1}//
E_{\Sigma X}))). \eea

Thus the $E_2$ term is concentrated in even degrees, and the spectral sequence
collapses.  Now we need to show its convergence.  Note that the computations
with the bar spectral sequences in the proof of Theorem
\ref{int:knfree} show that we have the following exact sequence of 
Hopf algebras: 
$$ E_X\otimes P_X \ra K(n)_*(\Pi _{i\in J}\und{\rm BP}{d_i})\ra
K(n)_*(\Omega ^{\infty }Y) \ra P_{\Sigma X} \ra K(n)_*.$$
On the other hand, we also have the following exact sequence from above:
$$E_2(1) \ra E_2(2) \ra E_2(3) \ra P_{\Sigma X} \ra K(n)_*.$$  Thus 
$E_2(3) $ is isomorphic to an associated graded object of 
$K(n)_*(\Omega ^{\infty }Y) $ so that the spectral sequence converges.  
Furthermore, as $E_2(3)$ is a polynomial algebra, there can be no 
nontrivial algebra extension,  which shows that $K(n)_*(\Omega ^{\infty }Y) $
is a polynomial algebra as desired.  
\hfill\qed\medbreak

Now we embed $K(n)_*(\Omega ^{\infty } C_f)$ into a more familiar object.

\proclaim{Proposition} \label{cok;mono}
Let $\{ g_i:C_f\ra \und{\rm BP}{e_i}|i\in J\} $ be a set with the property that
$K(n)_*(\Pi _ig_i):K(n)_*(\vee \Sigma ^{e_i}{\rm BP})\ra K(n)_*(C_f)$ is a monomorphism{\rm ,}
for example   a set of topological ${\rm BP}^*{\rm (BP)}$\/{\rm -}\/module
generators for ${\rm BP}^*(C_f)${\rm .}  Then $K(n)_*(\Pi _i\Omega ^{\infty }g_i):
K(n)_*(\Omega ^{\infty } C_f)\ra K(n)_*(\Pi _i\und{\rm BP}{e_i})$ is a monomorphism{\rm .}
\endproclaim

\demo{Proof}
By Propositions \ref{dll;quo} and  \ref{cof;pol},  we see that 
$K(n)_*(\Omega ^{\infty } \Sigma ^rC_f)$\break ($r\in Z^+$) forms a $Z^+$-indexed 
family of  free algebras with Dyer-Lashof length-like filtration.  Thus
one has the desired result.\enddemo

Thus we are ready to identify $K(n)_*(QX)$ as well as ${\rm BP}^*(QX)$.

\demo{Proof of Theorem {\rm \ref{knqx;alg}}} The first statement is obtained by combining Propositions \ref{res;debut} 
and \ref{cok;mono}.  The second follows from the first by Theorems 
\ref{land;gen} and \ref{qx;ten}.  The readers may object that 
$\Pi _j\und{\rm BP}{e_j}$ is not necessarily of finite type.  However, since this
space is torsion-free,  one can prove directly the properties from (v)$'$ to 
(xi) of Theorem 3.2.  Thus we can apply Theorem \ref{qx;ten} to the sequence
$QX\ra \Pi _i\und{\rm BP}{d_i} \ra \Pi _j\und{\rm BP}{e_j}$. 
\enddemo

\demo{Remark}  (i) Of course,  to obtain   part (i) of Theorem 3.2 for a fixed $n$,  it suffices
to assume that one has an exact sequence of the form
$0\ra \wt{K(n)_*}(X)\ra K(n)_*(\vee _i\Sigma ^{d_i}{\rm BP})\ra
K(n)_*(\vee _j\Sigma ^{e_j}{\rm BP})$. $\phantom{\sum^\int}$ \medbreak
\listitem{ii} The computation of $E_*(\und{\rm BP}{*})$ for complex oriented cohomology 
theories $E$ in \cite{RW} makes the evaluation of the map at the right end of
the exact sequence a completely algebraic process,  as explained in 
\cite{bjw}. It can also  be  reduced to evaluating $g_j$'s as cohomology 
operations on ${\rm BP}^*(CP^{\infty }\times \cdots \times CP^{\infty })$.  This can
be seen as follows.  According to \cite{kunst},  
$E_*(\und{\rm BP}{2*})$ is spanned by elements of the form $E_*(\alpha )(\beta )$
where $\alpha \in {\rm BP}^{2*}(CP^{\infty }\times \cdots \times CP^{\infty })\cong
[CP^{\infty }\times \cdots \times CP^{\infty },\und{\rm BP}{2*}]$,  $\beta \in
E_*(CP^{\infty }\times \cdots \times CP^{\infty })$.  But
$E_*(g_j)(E_*(\alpha )(\beta ))=E_*(g_j(\alpha ))(\beta )$,  where $g_j(\alpha )
$ denotes $g_j$ evaluated on $\alpha $.

\section{The destabilization functor for ${\rm BP}$-cohomology} \label{bp;des}

In this section we prove Lemma \ref{gen=pre} and complete the proof of Theorem
\ref{realmain}.  For this purpose, we start by studying the nature of the
skeleton filtration on ${\rm BP}^*(X)$ when it is well-generated.  First we improve
Lemma 4.3  of \cite{rwy}.

\proclaim{Lemma}
Let $X$ be a spectrum{\rm .}
Denote by $E^*_{*,*}(X)$ and $E^*_{*,*}({\rm sk}_mX)$ respectively the Atiyah\/{\rm -}\/Hirzebruch
spectral sequences $H^*(X,{\rm BP}^*)\Rightarrow {\rm BP}^*(X)$ and
$H^*({\rm sk}_mX,{\rm BP}^*)\Rightarrow {\rm BP}^*({\rm BP}_mX)${\rm .}
Then the natural map of the spectral sequences $E^r_{s,t}(X)\ra E^r_{s,t}({\rm sk}_mX)$ 
induced by the inclusion of the skeleton is an isomorphism for $s\leq m-r+1$
and a monomorphism for $m-r+2 \leq s \leq m${\rm . }
%\begin{enumerate}
%\item $E^r_{s,t}\ra E^r_{'s,t}$ is an isomorphism for $s\leq m-r+1$.
%\item $E^r_{s,t}\ra E^r_{'s,t}$ is a monomorphism for $m-r+2 \leq s \leq m$.
%\item $E'_r^{s,t}$ vanishes if $s>m$.
%\end{enumerate}
\endproclaim

{\it Proof}.  We proceed by induction on $r$.  The assertions are clearly true
when $r=2$.  Now suppose that they are true for $r$.
Consider the following commutative diagram.
$$
\begin{array}{ccccc}
\noalign{\vskip5pt}
 {E^r_{s-r,*}(X)} & \stackrel{d_r}{\lrar}&{E^r_{s,*}(X)}&\stackrel{d_r}{\lrar}&{E^r_{s+r,*}(X)}
 \\
\noalign{\vskip4pt}
\Big|&&\Big|&&\Big|\\
\noalign{\vskip-4pt}
\downarrow&&\downarrow&&\downarrow\\
 {E^r_{s-r,*}({\rm sk}_mX)}& \stackrel{d_r}{\srar} &{E^r_{s,*}({\rm sk}_mX)}&\stackrel{d_r}{\srar}
&{E^r_{s+r,*}({\rm sk}_mX)}
\end{array}
$$  
When $s\leq m-r$, by the induction hypothesis,  the right vertical arrow is a
monomorphism, and the other two vertical arrows are isomorphisms.  Thus after
taking the homology in the middle,  one sees that the induced map
$E^{r+1}_{s,t}(X)\ra E^{r+1}_{s,t}({\rm sk}_mX)$ is an isomorphism.  When $m-r+1\leq s\leq m$,
the induction hypothesis implies that the left vertical arrow is an isomorphism,
and the middle one is a monomorphism.  Therefore, by passing to the homology,
we see that the map $E^{r+1}_{s,t}(X)\ra E^{r+1}_{s,t}({\rm sk}_mX)$ is a monomorphism.  Thus
we conclude that the assertions hold for any $r$.  \hfill\qed

\proclaim{Lemma}
Let $X$ be a space or spectrum whose ${\rm BP}$ cohomology is well\/{\rm -}\/generated{\rm .}
Then an element $\alpha $ of ${\rm BP}^*(X)$ lies in the kernel of the map
${\rm BP}^*(X)\ra {\rm BP}^*({\rm sk}_nX)$ if and only if there exist  elements $a_i\in {\rm BP}^*${\rm ,}
$x_i\in {\rm BP}^{l_i}(X)$ such that $l_i>n$  for all $i$ and $
\alpha =\Sigma _i a_ix_i${\rm .}
\endproclaim

\demo{Proof} The ``if'' part is obvious.  So it suffices to prove the``only if''
part.  Let $\{ e_i |i\in \Lambda \}$ be a set of elements of ${\rm BP}^*(X)$ such that
$\{\rho _X(e_i) |i\in \Lambda \}$ is a basis for $M_X$.  Since ${\rm BP}^*(X)$ is 
well-generated,  $\{ e_i |i\in \Lambda \}$ generates ${\rm BP}^*(X)$.  
Thus if $\alpha $ is an element of ${\rm BP}^*(X)$, it can be written as
$$\alpha =\Sigma _i \beta _i e_i$$
with $\beta _i \in {\rm BP}^*$.  This implies that if it is in $\Ker ({\rm BP}^*(X)\ra {\rm BP}^*({\rm sk}_nX))$, then
$$\alpha '=\Sigma _{d_i\leq n} \beta _i e_i \in \Ker ({\rm BP}^*(X)\ra {\rm BP}^*({\rm sk}_nX))$$
since the other terms are obviously in the kernel.  Let 
\hbox{$d={\rm min}\{d_i|\beta _i\neq 0\}$.}\break
Thus in the Atiyah-Hirzebruch
spectral sequence for ${\rm BP}^*({\rm sk}_nX)$ we see\break that $\Sigma _{d_i=d}
\beta _i\otimes e_i$ is a boundary element and in $E'_{\infty }$ we have
the equality\break $\Sigma _{d_i=d}
\beta _i\otimes e_i=0$  and we have a nontrivial additive extension of the form
$$\Sigma _{d_i=d}
\beta _i e_i =-\Sigma _{d<d_i\leq n}\beta _i e_i .$$
On the other hand,  as  we have seen previously, the map $E^{\infty }_{s,t}(X)
\ra E^{\infty }_{s,t}({\rm sk}_nX)$ is injective for $s\leq n$.  Thus by the 
naturalness of the Atiyah-Hirzebruch spectral sequence,  one can conclude that
in $E^{\infty }_{s,t}(X)$ there is the same type of extension problem, and
$$\Sigma _{d_i=d}
\beta _i e_i =-\Sigma _{d<d_i\leq n}\beta _i e_i 
+\Sigma _{d_i>n} \gamma _i e_i.$$  Thus we have 
$\alpha =\Sigma _{d_i>n} \gamma _i e_i+\Sigma _{d_i>n} \beta _i e_i$ as desired.  \hfill\qed
\medbreak

As a straightforward consequence, we have the following:

\proclaim{Proposition} \label{filt:good}
Let $X$, $Y$ be spectra or spaces and $f:X\ra Y$ a map such that ${\rm BP}^*(X)$ and
${\rm BP}^*(Y)$ are well\/{\rm -}\/generated and such  that ${\rm BP}^*(f)$ is surjective{\rm .}  Then the 
skeletal filtration on ${\rm BP}^*(X)$ agrees with the quotient filtration induced from
the skeletal filtration on ${\rm BP}^*(Y)${\rm .}
\endproclaim

\demo{Proof of  Lemma {\rm \ref{gen=pre}}} The case of a space follows from the case
of a spectrum.  When $X$ is a spectrum,  it suffices to consider a familly of
generators $X  \stackrel{\vee _if_i}{\srar} \vee _i\Sigma ^{d_i}{\rm BP}$,  and apply the lemma above.
\enddemo

Now we are ready to complete the proof of Theorems \ref{realmain} and \ref{evencell}.  As we need to deal with a topology on ${\rm BP}^*(-)$ that is different from
the skeletal topology,  we recall:

\numbereddemo{Definition}
The finite-subcomplex topology on ${\rm BP}^*(X)$ (whether $X$ is a space or a spectrum)
is the topology in which  the system of neighbourhoods of $ 0$ is the set of 
$\Ker ({\rm BP}^*(X)\ra {\rm BP}^*(X_{\alpha }))$ where $X_{\alpha }$ runs through all finite
subcomplexes of $X$.
\enddemo

\demo{Remark}
This topology is often called pro-finite topology in the litterature (e.g.\
\cite{ad}, \cite{boa}).  We prefer to rename it because the term pro-finite topology means
something else for algebraists (the topology in which the neighbourhoods of 0
are the subgroups of finite index),  and it is not too absurd to consider
the pro-finite topology in this sense here.

\demo{Proof of Theorem \ref{realmain}}  
According to Theorem \ref{knqx;alg},  we have the exact sequence of ${\rm BP}^*{\rm (BP)}$-modules 
$$0\la \widetilde{\rm BP}^*(X)\la {\rm BP}^*(Y_1) \la {\rm BP}^*(Y_2)$$ and a coexact sequence of
augmented ${\rm BP}^*$-algebras 
$$ {\rm BP}^* \la  {\rm BP}^*(QX) \la {\rm BP}^*(\Omega ^{\infty }Y_1)\la {\rm BP}^*(\Omega ^{\infty }Y_2),
$$ where $Y_1$ and $Y_2$ are wedges of suspensions of ${\rm BP}$.
On the other hand, $X$ being of finite type, one can take $Y_1$ to be of finite
type as well,  which forces $\Omega ^{\infty }Y_1$  to be of finite
type also.  Thus one sees that ${\rm BP}^*(Y_1)$ is free and ${\rm BP}^*(\Omega ^{\infty }Y_1)$
can be identified with ${\cal D}({\rm BP}^*(Y_1))$.  Unfortunately $Y_2$ is not 
necessarily of finite type.  If $Y_2\cong \vee _j \Sigma ^{e_j}{\rm BP}$,
${\rm BP}^*(Y_2)$ is the completion with respect to the finite subcomplex topology of
$\mathbold{\oplus}_j \Sigma ^{e_j}{\rm BP}^*{\rm (BP)}$.  It contains the completion with respect to
the skeletal topology of $\oplus _j \Sigma ^{e_j}{\rm BP}^*{\rm (BP)}$, which we call $M$.
Thus $M$ is free, dense in ${\rm BP}^*(Y_2)$ with respect to the finite subcomplex
topology,  and complete with respect to the skeletal topology.  Since these two
topologies are natural,   $\Im(M\ra {\rm BP}^*(Y_1))$  is dense in $\Im({\rm BP}^*(Y_2)
\ra {\rm BP}^*(Y_1))$ with respect to the finite subcomplex topology,  and complete with
respect to the skeletal topology.  But these two topologies agree on ${\rm BP}^*(Y_1)$,
which implies that $\Im(M\ra {\rm BP}^*(Y_1))$ coincides with $\Im({\rm BP}^*(Y_2)
\ra {\rm BP}^*(Y_1))$.  Thus we can replace the first exact sequence by the following 
one:
$$0 \la \widetilde{\rm BP}^*(X)\la {\rm BP}^*(Y_1) \la M.$$
Here the two terms on the right are free.  
  On the other hand,  Lemma \ref{gen=pre}
shows that the skeletal filtration on $\widetilde{\rm BP}^*(X)$ agrees with the 
quotient filtration induced from that of ${\rm BP}^*(Y_1)$.  Thus this is really an
exact sequence in
 ${\cal M}'_{\rm BP}$.
Using  the definition of ${\cal D}$  we can identify
${\rm BP}^*(\Omega ^{\infty }Y_1)$ with ${\cal D}{\rm BP}^*(Y_1)$.  Furthermore,
by arguments similar to the one above,  we can 
replace ${\rm BP}^*(\Omega ^{\infty }Y_2)$ by ${\cal D}(M)$ in the coexact sequence 
above.  Again  Proposition 
\ref{filt:good} shows that ${\rm BP}^*(QX)$ has the quotient filtration and that 
the sequence is coexact in ${\cal K}' _{_0{\rm BP}}$.  
Since ${\cal D}$ has to be right exact,  we obtain the desired result.   \enddemo

We can also generalize the result in \cite{bpqs}.

\demo{Proof of Theorem {\rm \ref{evencell}}} The case in which $i$ is even is 
essentially
already treated in \cite[Cor.\  to Th.\ 7.3]{bpqs}.  However, there is one mistake in the proof 
in \cite{bpqs}. It was implicitely assumed there that ${\rm BP}$-cohomology of a
wedge of suspensions of ${\rm BP}$ is  a free module over ${\rm BP}^*{\rm (BP)}$,  which is false
unless  we take into account the topology (thus replacing $D$ defined 
there by ${\cal D}$ defined here) or if we assume that ${\rm BP}^*(X)$ is finitely
generated as a ${\rm BP}^*{\rm (BP)}$-module so that we do not have to worry about the
topology.  

Now
let $i=2j-1$.  Then as was shown in $\cite{bcat}$ one has the inclusion
${\rm Ind} K(n)_*(\und{X}{2j-1})\subset K(n)_*(\und{X}{2j})$.  Since 
$PK(n)_*(\und{X}{2j-1})\cong {\rm Ind} K(n)_*(\und{X}{2j-1})$,  in cohomology we 
have that the map $K(n)^*(\und{X}{2j})\ra {\rm Ind} 
K(n)^*(\und{X}{2j-1})$ is an epimorphism.
Since ${\rm BP}^*(\und{X}{2j})\bpg K(n)^*$ surjects to $K(n)^*(\und{X}{2j})$,
one sees that ${\rm BP}^*(\und{X}{2j-1})\bpg K(n)^*$ surjects to 
$K(n)^*(\und{X}{2j-1})$.  However,  the arguments in \cite{bpqs} show that we
have all the exact sequences of Hopf algebras needed in Morava $K$-theories,
so that Theorem \ref{land;gen} implies the desired result. \phantom{romeo} \enddemo

\appendix{}{Modifications for prime $2$}

In this appendix we treat the case $p=2$.  First of all,  we don't have to worry
about the possible noncommutativity of Morava $K$-theory,  since in cohomology,
all spaces  dealt with satisfy ${\rm BP}^*(X)\bpg K(n)\cong K(n)^*(X)$,  so that
the cup product is commutative.  In homology,  all $H$-spaces  dealt with have
$H$-maps from each one to another space whose Morava $K$-homology is known to be 
commutative; these maps induce  monomorphism in Morava $K$-homology,  so 
that the Morava $K$-homology of these $H$-spaces is commutative.  There remain two
sources of   problems.  First of all,  the square of odd degree elements in
commutative graded $Z/2$-algebras is not necessarily zero,  which requires us to
revise the content of Sections \ref{fil;dll} and \ref{kn;poly}.  Another thing is that 
the Adem relations,  May's formula,  and Nishida relations do not exactly look
the way they do when $p$ is odd,  which makes us  modify the arguments in 
Section~\ref{thom;im} a little bit.  Now we list what changes.

In Section \ref{prel}, we first  replace Definition \ref{def;11}.

\numbereddemo{Definition} 
$I=(s_1,\dots ,s_k)$ is called admissible if for $s_j\leq 2s_{j+1}${\rm ,}  the excess{\rm ,}  the  degree{\rm ,}  and
the length of $I$ are defined by $d(I)=\Sigma _{j=1}^ks_j${\rm ,}  $l(I)=k$ and $e(I)=s_1-\Sigma
_{j=2}^ks_j${\rm .} \enddemo 

With this modification,  Theorem \ref{th;12} holds as stated,  except that now
the relevant reference is \cite{AK}.  The formulae in Theorem \ref{mayf} which
do not involve $\beta $ hold by replacing $P^i$ with $Sq^i$.  In particular we 
have $\beta Q^{2s}=Q^{2s-1}$.  Observe also that the formulae in Theorem 
\ref{mayf} which involve $\beta $ hold as well by replacing $Q^i$ with $Q^{2i}$
and $P^i$ with $Sq^{2i}$.  
Finally Proposition 
\ref{prop;14} is true modulo the algebra extension $x^2=y$.

In Section \ref{fil;dll}, everything remains valid if we replace ``free
commutative algebra" with ``tensor product of a polynomial algebra concentrated
in even degrees and an exterior algebra generated by odd degree elements".  The
results in this section now can be used in Section \ref{kn;poly}.  However,  we
need  a variant of Proposition \ref{dll;main} that can be used in the proof of
Lemma \ref{tho;lb}.  For this purpose,  we change conditions (i) and (iv) of the
Definition \ref{filt;dll} as follows:
\medbreak
\listitem{i$'$} Each  $A_{i,*}$ is a polynomial algebra,  and  there exists $\varepsilon $,
such that $A_{i,*}$ is generated by even degree elements if and only if 
$i$ is congruent to $\varepsilon$ mod $2$,  and such that $A_{i,*}$ is generated by 
odd degree elements otherwise. \smallbreak 
\listitem{iv$'$} $\sigma $ induces an isomorphism ${\rm Ind}(A_{2i+\varepsilon })\ra
{\rm Ind}(A_{2i+\varepsilon +1})$ and a monomorphism ${\rm Ind}(A_{2i+\varepsilon -1})\ra
A_{2i+\varepsilon }$.  
\medbreak

Now with this definition for the algebras with Dyer-Lashof length-like 
filtration,  a variant of Proposition \ref{dll;main} holds by requiring $B$ to 
satisfy the following conditions:
\medbreak
\listitem{i} Each $B_{2i+\varepsilon }$ is a polynomial algebra concentrated
in even degrees. \smallbreak
\listitem{ii} Each $B_{2i+\varepsilon -1 }$ is an exterior algebra generated by
odd degree elements.
\medbreak
  The proof is similar to the odd prime case.  

In Section \ref{thom;im}, the  following modifications are required.  Wherever we
consider an operation $Q^I$   containing at least one 
Bockstein for each odd prime,  we consider an operation $Q^I$,  $I=(s_1,\dots s_l)$,
with at least one $s_j$ being odd.  Then the proof of  Proposition 
\ref{bpg;det} can be proved in a similar way as in the  odd prime case, taking into
account the observation we made after the modifications on Theorem \ref{mayf}.
Finally we prove a weakened version of Lemma \ref{tho;lb}.  
First we~show:

\proclaim{Lemma} \label{new52}
Let $\{ f_i:X\ra \und{{\rm BP}}{d_i}|i\in I\} $ be a set 
that reduces to a $Z/2$-basis for $M_X\subset \wt{H}^*(X)${\rm .}
Then one has
$$\Im (\otimes _{i\in I}H^*(\Omega ^{\infty }f_i):
H^*(\Pi _{i\in I}\und{{\rm BP}}{d_i})\ra H^*(QX))=C$$
where $C$ is as in   Proposition {\rm \ref{bpg;det}.}
\endproclaim

\demo{Proof}
We proceed as in the odd prime case to take a subalgebra $T_{\Sigma ^nX}$
of $H_*(Q\Sigma ^nX)$ similarly.  If we give an increasing filtration on
 $T_{\Sigma ^nX}$ by defining $F_j(T_{\Sigma ^nX})$ to be the subalgebra 
generated by the image of Dyer-Lashof operations on the elements in 
$H_j({\Sigma ^nX})$,  and that on $H_*(\Pi _{i\in I}\und{{\rm BP}}{d_i+n})$  by
defining $$F_j(\Pi _{i\in I}\und{{\rm BP}}{d_i+n})\cong H_j(\Pi _{i\in I,d_i+n\leq j}
\und{{\rm BP}}{d_i+n}),
$$ then $H_*(\Pi _if_i)$ respects this filtration.  
However,  $F_{l}//F_{l-1}$ of the former is just $(A_{\Sigma ^nX})_{l}\otimes
T_{\Sigma ^nS^{l-n}}$ whereas that of the latter is $(A_{\Sigma ^nX})_{l}\otimes
H_*(\und{{\rm BP}}{l})$.  We know that the former injects to the latter either by our
modified version of Proposition \ref{dll;main} or by \cite{W}.  Thus we see that
$T_X$ injects to $H_*(\Pi _{i\in I}\und{{\rm BP}}{d_i})$.  The rest of the proof does
not require modification.  \enddemo

Fewer $f_i$'s will suffice,  at least when $X$ satisfies the hypotheses of
Theorem \ref{realmain}.  We will come back to this point later.

In Section \ref{kn;poly},  we first  have to weaken Theorem \ref{kn;mon}.  
Namely we have to take $\{ f_i:X\ra \und{{\rm BP}}{d_i}|i\in I\} $ be a set
that reduces to a $Z/2$-basis for $M_X\subset \wt{H}^*(X)$.  Then using Lemma
\ref{new52} instead of Lemma \ref{tho;lb},  one can prove the theorem in a 
similar way.  Next,  throughout the section,  ``polynomial algebra" should be 
replaced with ``polynomial algebra concentrated in even degrees" and 
``exterior algebra" should be 
replaced with ``exterior algebra generated by odd degree elements".  Then the rest of 
the section becomes true.  Note that using Proposition \ref{kn;mono},  one sees
that Theorem \ref{kn;mon} holds without changing the family 
$\{ f_i:X\ra \und{{\rm BP}}{d_i}|i\in I\} $ under the hypotheses of Theorem \ref{realmain}.
Thus under these assumptions,
one can prove Lemma \ref{tho;lb} for the original family 
$\{ f_i:X\ra \und{{\rm BP}}{d_i}|i\in I\} $ in the statement.

Throughout Section \ref{kn;poly},  ``free commutative
algebra" should be
replaced with ``polynomial algebra concentrated in even degrees tensored with
exterior algebra generated by odd degree elements",  and ``cofree cocommutative
coalgebra" with ``divided power coalgebra concentrated in even degrees tensored
 with exterior coalgebra generated by odd degree elements".  Some of the spectral 
sequences may have possible nontrivial algebra  extension problems due to 
the fact that the squares of odd degree elements are not automatically zero.
However,  using the naturality arguments and comparison with appropriate 
spectral sequences,  one can always show that the statements in this section
remain true after the modification mentioned above.

All the rest of the article remains true as stated.

\bye

\item"[{E}]"   
D.\ B.\ A. Epstein, Complex hyperbolic geometry, {\it London Math.\
Soc. Lecture Notes\/} {\bf 111} (1987), 93--111.

\item"[{FZ}]"
E.\ Falbel and V.\ Zocco, A Poincar\'e's fundamental polyhedron
theorem for complex hyperbolic manifolds, {\it J.\ Reine Angew.\
Math\/}.\
{\bf 516} (1999), 133--158.

\item"[{G}]"
W.\ Goldman, {\it Complex Hyperbolic Geometry\/}, Oxford Univ.\ 
Press, New York, 1999.

\item"{[GKL]}"
W.\ Goldman, M.\ Kapovich, and B.\ Leeb, Complex hyperbolic surfaces
homotopy equivalent to a Riemann surface, {\it Comm.\ Geom.\ Anal\/}.,
to appear.

\item"{[GMT]}"
D.\ Gabai, R.\ Meyerhoff,and N.\ Thurston, Homotopy hyperbolic
manifolds are hyperbolic, {\it Ann.\ of Math\/}., to appear.

\item"{GP]}"
W.\ Goldman, and J.\ Parker, Complex hyperbolic ideal triangle groups,
{\it J.\ Reine Angew.\ Math\/}.\ {\bf 425} (1992), 71--86.

\item"{[GuP]}"
N.\ Gusevskii and J. Parker, Complex hyperbolic quasifuchsian
surfaces, preprint, 1999.

\item"{[K]}"
F.\ Klein, Neue Beitr\"age zur Riemannschen functionentheorie, {\it 
Math.\ Ann\/}.\ {\bf 21} (1883), 141--218.

\item"{[KR]}"
A.\ Kornayi and H.\ M.\ Riemann, Quasiconformal mappings on the
Heisenberg group, {\it Invent.\ Math\/}.\ {\bf 80} (1985), 309--338.

\item"{[KeR]}"
B.\ Kernighan and D.\ Ritchie, {\it The C Programming Language\/},
Prentice-Hall, 1978.

\item"{[I]}"
{\it IEEE standard for binary floating-point arithmetic\/}, 
Inst.\ of Electrical and Electronics Engineers, July 26, 1985.

\item"[Mac]"
A.\ M.\ Macbeath, Packings, free products and residually finite
groups, {\it \ Proc.\ Cambridge Philos.\ Soc\/}.\ {\bf 59} (1963),
555--558.

\item"[Mas]"
B.\ Maskit, Construction of Kleinian groups, {\it Proc.\ Conf.\
Complex Analysis\/} (Minneapolis, 1964), 281--296, Springer-Verlag, New
York, 1965.

\item"[Sa]"
H.\ Sandler, Trace equivalence in SU$(2,1)$, {\it Geom.\ Dedicata\/}
{\bf 69} (1998), 317--327.

\item
"[{S1}]"
R.\ Schwartz, Dented tori, 1997, electronic program-document,
written in TcL and C, available upon request to res\@math.umd.edu.

\item
"[{S2}]"
R.\ Schwartz, Degenerating the complex hyperbolic ideal triangle groups,
{\it Acta Math\/}., to appear.

\item
"[{S3}]"
R.\ Schwartz, Circle quotients and string art, preprint, 2000.

\item
"[S4]"
R.\ Schwartz, Real hyperbolic on the outside, complex hyperbolic on the inside,
preprint, 2000.

\item"[T]"
W.\ Thurston, On the geometry and dynamics of diffeomorphisms
of surfaces, {\it Bull.\ A.M.S\/}.\ {\bf 19} (1988), 417--431.

\item"[Tol]"
D.\ Toledo, Representations of surface groups in complex hyperbolic
space, {\it J.\ Differential Geom\/}.\ {\bf 29} (1989), 125--133.

\item"[W]"
S.\ Wolfram, Mathematica: A system for doing mathematics by
computer, 1995.

\item
"[{W-G}]"
J.\ Wyss-Gallifent, Complex hyperbolic triangle groups, Ph.D.\
thesis, Univ. of Maryland, 2000.

\endroster

\bigskip

\centerline{(Received March 16, 1998)}
\centerline{(Revised August 17, 2000)}

\endreferences

\begin{references}

 



\bibitem{ad} \name{J. F. Adams},  
A variant of E. H. Brown's representability
theorem, {\it Topology\/} {\bf 10} (1971), 185--198.

\bibitem{AK} \name{S.\ Araki} and \name{T.\ Kudo},   Topology of $H_n$-spaces and
$H$-squaring operations, {\it Mem.\  Fac.\ Sci.\ Kyusyu
Univ.\ Ser.\ A\/} (1956),  85--120.

\bibitem{be} \name{M.\  G.\ Barratt} and \name{P.\ J.\ Eccles},  
 $\Gamma ^+$-structures-III.\ 
 The stable structure of $\Omega ^{\infty }\Sigma ^{\infty }A$, {\it Topology\/}
{\bf 13} (1974), 199--207.

\bibitem{bcm} \name{M. Bendersky,  E. B. Curtis,} and \name{H. R. Miller}, 
The  unstable Adams spectral sequence for generalized homology,  {\it
Topology\/} {\bf 17} (1978), 229--248.

\bibitem{boa} \name{J. M. Boardman},
Stable operations in generalized cohomology,   
in {\it The Handbook of Algebraic Topology},  
North-Holland, Amsterdam, 585--686, 1995.

\bibitem{bjw} \name{J. M. Boardman,  D. C. Johnson,} and \name{W. S.
Wilson}, Unstable operations in generalized cohomology,   
{\it The Handbook of Algebraic Topology}, 687--828,  North-Holland, Amsterdam,  1995.

\bibitem{BKf} \name{A. K. Bousfield},   On
$p$-adic  $\lambda$-rings and the $K$-theory of 
$H$-spaces,  {\it Math. Z\/}.\  {\bf 223} (1996), 483--519.

\bibitem{BK} \bibline,   On $\lambda$-rings and the $K$-theory of 
infinite loop spaces,  $K$-{\it Theory\/} {\bf 10} (1996), 1--30.

\bibitem{dl} \name{E.\ Dyer}  and \name{R.\ K.\ Lashof},  Homology of iterated loop 
spaces, {\it  Amer.\ J.\  Math\/}.\  {\bf 84} (1962),  35--88.

\bibitem{Ho} \name{L. Hodgkin},  The $K$-theory of some well-known
spaces.  I. $QS ^0$, {\it Topology\/} {\bf 11}  (1972),  371--375.

\bibitem{Ho2} \bibline,   Dyer-Lashof operations in $K$-theory,  in {\it
New  Developments in Topology\/} ({\it Proc.\ Sympos.\ Algebraic
Topology, Oxford\/}, 1972), 27--32,
{\it London Math.\ Soc.\ Lecture Note Ser\/}., No.\ 11, Cambridge Univ.\ Press, London, 1974.


\bibitem{hs} \name{L.\ Hodgkin} and \name{V.\ P.\ A.\ Snaith},   The $K$-theory of some
more well-known spaces,  {\it Illinois J.\ Math\/}.\ {\bf 22} (1978),
270--278.

\bibitem{hh}  \name{M.\ J.\ Hopkins} and  \name{J.\ R.\ Hunton},   On the structure of
spaces representing a Landweber exact cohomology theory,  {\it
Topology\/} {\bf 34} (1995), 29--36.

\bibitem{HKR}
\name{M.\ J.\ Hopkins,  N.\ J.\ Kuhn},  and \name{D.\ C.\ Ravenel},   Generalized group 
characters and complex
oriented cohomology theories, {\it J.A.M.S.} {\bf 13} (2000), 553--594.


\bibitem{Hu1} \name{J. R. Hunton},  Morava $K$-theory of wreath
products, {\it Math.\ Proc.\ Cambridge Philos. Soc\/}.\ {\bf 107}
(1990),   309--318.

\bibitem{huad} \bibline, Detruncating Morava $K$-theory, in
{\it Adams Memorial Symposium on Algebraic Topology\/} {\bf 2}
(Manchester, 1990), {\it
 London Math.\
Society Lect.\ Note Ser\/}.\  {\bf 176},  Cambridge Univ.\  Press,  London, 1992.


\bibitem{jrh;gen} \bibline,  The complex oriented cohomology of
extended powers,  {\it Ann.\  Inst.\ Fourier\/} {\bf 48} (1998), 517--534.


\bibitem{In}  \name{K. Inoue},  The Brown-Peterson cohomology of
$B{\rm SO}(6)$, {\it J.\ Math.\  Kyoto Univ\/}.\ {\bf 32} (1992), 655--666.


\bibitem{jw} \name{D.\ C.\ Johnson} and \name{W.\ S.\ Wilson},  ${BP}$ operations and
Morava's extraordinary $K$-theories,  {\it Math.\ Z\/}.\ {\bf 144}
(1975), 55--75.

\bibitem{K1} \name{T. Kashiwabara},  Mod $p$ $K$-theory of
$\Omega^{\infty}
\Sigma^{\infty }X$ revisited,  {\it Math.\ Proc.\ Cambridge Philos.\
Soc\/}.\
 {\bf 114}
(1993), 219--221.

\bibitem{kunst} \bibline,   Hopf rings and unstable
operations,   {\it J.\  Pure  Appl.\ Algebra\/}  {\bf 94} (1994),  183--193.

 \bibitem{k2} \bibline, $K$(2)-homology of some infinite loop spaces,
{\it Math. Z\/}.\  {\bf 218} (1995), 503--518.

\bibitem{bpqs} \bibline,   Brown-Peterson cohomology
of
$\Omega^\infty\Sigma^\infty S^{2n}$,
{\it Quart.\  J.\  Math\/}.\ {\bf 49} (1998), 345--362.

\bibitem{Kprep} \bibline, Unstable Hopf algebras, destabilization
and cohomology of infinite loop spaces, in preparation.


\bibitem{bcat}  \name{T.\ Kashiwabara,  N.\ P.\ Strickland},  and 
\name{P.\ R.\ Turner}, 
The Morava $K$-theory Hopf ring for ${BP}$,
in {\it Algebraic Topology\/}: {\it New Trends in Localization and
Periodicity\/},
{\it Progr.\ Math\/}.\ {\bf 136}, 
 Birkh\"{a}user, Basel, 1996.


\bibitem{KW} \name{T.\ Kashiwabara} and \name{W.\ S.\ Wilson},  The Morava
$K$-theory and Brown-Peterson cohomology of spaces related ${BP}$, 
preprint.



\bibitem{LZ}  \name{J.\ Lannes} and \name{S.\ Zarati},   Sur les foncteurs 
d\'{e}riv\'{e}s de la d\'{e}stabilisation,  {\it Math.\ Z\/}.\
{\bf 194} (1987), 25--59.


\bibitem{geo}   \name{J.\ P.\ May}, {\it  The Geometry of Iterated Loop
Spaces\/},  {\it Lect.\ Notes in Math\/}.\  {\bf 271}, Springer-Verlag,
New York, 1972.

\bibitem{May} \bibline,   The homology of $E_{\infty}$ spaces, in
    {\it The Homology of Iterated Loop Spaces\/},
  {\it Lect.\  Notes in Math\/}.\  {\bf 533},  1--68, Springer-Verlag,
New York, 1976.

\bibitem{mc}  \name{J. E. McClure},  In mod $p$ $K$-theory of $QX$, in 
{\it $H_{\infty}$ Ring
Spectra and their Applications},   {\it Lect.\  Notes in Math\/}.\ {\bf
1176}, 291--383,   
Springer-Verlag, New York,  1986.

\bibitem{ms} \name{H.\ R.\ Miller} and \name{V.\ P.\ A.\ Snaith},   On 
$K_*(Q{\bf  R}P^n;P{\bf Z}/2)$,  in {\it Current Trends in
Algebraic Topology, Part \/} I, {\it CMS Conf.\ Proc\/}.\ {\bf 2}
(1982),  233--243.

\bibitem{Q} \name{D. Quillen},  Elementary proofs of some results of
cobordism theory using Steenrod operations, {\it Adv.\  in Math\/}.\ {\bf
7}  (1971), 29--56.


\bibitem{RW}
 \name{D.\ C.\ Ravenel} and \name{W.\ S.\ Wilson},   The Hopf ring for complex cobordism,
{\it J.\  Pure Appl.\  Algebra\/} {\bf 9}  (1977),  241--280.

\bibitem{rwy} \name{D.\ C. Ravenel,  W.\ S.\ Wilson},  and \name{N.\
Yagita},  
 Brown-Peterson cohomology from Morava $K$-theory, 
$K$-{\it Theory\/} {\bf 15}
(1998), 147--199.

\bibitem{Lsm} \name{L. Smith}, On the construction of the Eilenberg-Moore
spectral sequence, {\it Bull.\ A.\ M.\ S\/}.\ {\bf 75} (1969), 873--878.

\bibitem{sn}  \name{V. P. A. Snaith},   A stable decomposition of 
$\Omega^n S^nX$, {\it J. London Math.\  Soc\/}.\  {\bf 7} (1974), 
577--583.


\bibitem{SE}  \name{N. E. Steenrod},  {\it Cohomology
Operations\/},
{\it Ann.\ of Math.\ Studies\/}, No.\ 50, Princeton Univ.\  Press,
Princeton, NJ, 1962.

\bibitem{neil} \name{N. P. Strickland}, unpublished.

\bibitem{S} \bibline,  Morava $E$-theory of symmetric groups,
  {\it Topology\/} {\bf 37} (1998), 757--779; Correction: {\it ibid\/}.\
  {\bf 38} (1999), 931.

\bibitem{ST} \name{N.\ P.\ Strickland} and \name{P.\ R.\ Turner},   Rational Morava
$E$-theory of $DS^0$,   {\it Topology\/} {\bf 36} (1997),  137--151.


\bibitem{Ta}  \name{D. Tamaki},
 A dual Rothenberg-Steenrod spectral sequence,  {\it Topology\/} 
 {\bf 33} (1994),
631--662.

\bibitem{W} \name{W. S. Wilson},    The $\Omega$-spectrum for  
Brown-Peterson cohomology  I,  {\it Comment. Math. Helv\/}.\  {\bf 48}
(1973),  45--55; corrigendum: {\it ibid\/}.\ {\bf 48} (1993), 194.

\bibitem{wbook} \bibline, {\it Brown-Peterson Homology\/}:  
{\it An Introduction and 
Sampler\/}, {\it CBMS Reg.\ Conf.\ Ser.\ Math}.\
{\bf 48},  Conf.\ Board of Math.\ Sciences, Washington, D.C., 1982.

\bibitem{w2} \bibline,    Brown-Peterson cohomology from
Morava $K$-theory II,   $K$-{\it Theory\/} {\bf 17} (1999),  95--101.
\end{references}
\bye